\documentclass[a4paper,10pt]{scrartcl}
\usepackage[english]{babel}
\usepackage[utf8]{inputenc}
\usepackage{enumerate}
\usepackage{geometry}
\usepackage{amsfonts,amsmath,amssymb,amsopn}
\usepackage{amsthm}
\usepackage{mathtools}
\usepackage{stmaryrd}
\usepackage{nicefrac}
\usepackage{exscale}

\usepackage[numbers,square]{natbib}
\usepackage{url} 
\usepackage[colorlinks=true,pdfpagelabels,unicode]{hyperref}
\hypersetup{linkcolor=blue, urlcolor=blue, citecolor=blue} 

\allowdisplaybreaks 

\theoremstyle{plain}

\newtheoremstyle{theo}
	{3pt} 
	{3pt} 
	{\itshape} 
	{} 
		{\bfseries} 
	{\\} 
	{ } 
	{\thmname{#1}\thmnumber{ #2.}\thmnote{ - #3}} 
\theoremstyle{theo}

\newtheorem{definition}{Definition}[section]
\newtheorem{lemma}[definition]{Lemma}
\newtheorem{theorem}[definition]{Theorem}
\newtheorem{corollary}[definition]{Corollary}
\newtheorem{remark}[definition]{Remark}

\newenvironment{bew}{\begin{proof}[\bfseries Proof:]}{\end{proof}}

\DeclareMathOperator{\bomega}{\overline{\Omega}}
\DeclareMathOperator{\romega}{\partial\Omega}

\DeclareMathOperator{\intd}{d\!}

\DeclareMathOperator{\dive}{\nabla\cdot}
\DeclareMathOperator{\wto}{\rightharpoonup}

\newcommand{\epsi}{\varepsilon}

\usepackage[usenames,dvipsnames,svgnames,table]{xcolor}

\newcommand{\Tme}{T_{max,\;\!\epsi}}
\newcommand{\GNI}{Gagliardo--Nirenberg inequality}

\newcommand{\into}[1]{\int_0^{#1}\!}

\newcommand{\intoT}{\into{T}}

\newcommand{\intomega}{\int_{\Omega}\!} 

\newcommand{\intoTomega}{\intoT\!\intomega}
\newcommand{\intinfomega}{\int_0^\infty\!\!\intomega}

\newcommand{\Lo}[1][1]{L^{#1}(\Omega)} 
\newcommand{\W}[1][1,2]{W^{#1}(\Omega)}

\newcommand{\LSp}[2]{L^{#1\;\!}\!\left(#2\right)} 
\newcommand{\LSploc}[2]{L_{loc}^{#1}\!\left(#2\right)} 
\newcommand{\WSp}[2]{W^{#1}\!\left(#2\right)}
\newcommand{\CSp}[2]{C^{#1}\!\left(#2\right)}

\newcommand{\R}{\mathbb{R}}
\newcommand{\N}{\mathbb{N}}

\newcommand{\nfrac}[2]{{\nicefrac{#1}{#2}}}

\author{Tobias Black\thanks{Institut f\"ur Mathematik, Universit\"at Paderborn, Warburger Str. 100, 33098 Paderborn, Germany; email: \mbox{tblack@math.upb.de}}}
\title{Global very weak solutions to a chemotaxis-fluid system with nonlinear diffusion}
\date{}

\begin{document}
\maketitle
\begin{abstract}
\noindent
{\textbf{Abstract:} We consider the chemotaxis-fluid system
\begin{align}\label{star}\tag{$\diamondsuit$}
\left\{
\begin{array}{r@{\,}c@{\,}c@{\ }l@{\quad}l@{\quad}l@{\,}c}
n_{t}&+&u\cdot\!\nabla n&=\Delta n^m-\nabla\!\cdot(n\nabla c),\ &x\in\Omega,& t>0,\\
c_{t}&+&u\cdot\!\nabla c&=\Delta c-c+n,\ &x\in\Omega,& t>0,\\
u_{t}&+&(u\cdot\nabla)u&=\Delta u+\nabla P+n\nabla\phi,\ &x\in\Omega,& t>0,\\
&&\nabla\cdot u&=0,\ &x\in\Omega,& t>0,
\end{array}\right.
\end{align}
in a bounded domain $\Omega\subset\mathbb{R}^3$ with smooth boundary and $m>1$. Assuming $m>\frac{4}{3}$ and sufficiently regular nonnegative initial data, we ensure the existence of global solutions to the no-flux-Dirichlet boundary value problem for \eqref{star} under a suitable notion of very weak solvability, which in different variations has been utilized in the literature before. Comparing this with known results for the fluid-free setting of \eqref{star} the condition appears to be optimal with respect to global existence. In case of the stronger assumption $m>\frac{5}{3}$ we moreover establish the existence of at least one global weak solution in the standard sense.

In our analysis we investigate a functional of the form $\int_{\Omega}\! n^{m-1}+\int_{\Omega}\! c^2$ to obtain a spatio-temporal $L^2$ estimate on $\nabla n^{m-1}$, which will be the starting point in deriving a series of compactness properties for a suitably regularized version of \eqref{star}. As the regularity information obtainable from these compactness results vary depending on the size of $m$, we will find that taking $m>\frac{5}{3}$ will yield sufficient regularity to pass to the limit in the integrals appearing in the weak formulation, while for $m>\frac{4}{3}$ we have to rely on milder regularity requirements making only very weak solutions attainable.
\noindent
}\\[0.1cm]


\end{abstract}


\newpage
\section{Introduction}\label{sec1:intro}
The influence of chemotaxis, that is the biased movement of cells in the direction of chemical concentration gradients, on the evolution of cell populations has been one of the focal points of mathematical biology since the introduction of the acclaimed model
\begin{align}\label{KS}
n_t=\nabla\cdot(D(n)\nabla n- S(n,c)\nabla c)&& c_t=\Delta c-c+n
\end{align}
by Keller and Segel (\cite{KS70}), where $n(x,t)$ and $c(x,t)$ denote the density of the cell population and the concentration of the attracting chemical substance, respectively, at place $x$ and time $t$. The system is able to describe the spontaneous aggregation process of bacteria, which can be observed for populations of e.g. \emph{Dictyostelium discoideum}, and hence piqued the interest of many mathematicians. This fascinating behavior already emerges for the prototypical choices $D(n)\equiv1$ and $S(n,c)\equiv n$ if either the initial mass of cells $\intomega n_0$ is large enough (\cite{HVblow97}), or for certain initial data of arbitrary initial mass in dimensions $N\geq3$ (\cite{Win13pure}).

Biologically a stronger nonlinear diffusion, e.g. a porous medium type $D(n)\simeq mn^{m-1}$, seems appropriate as cells cannot be compressed to a single point and hence densely packed cells suffer a larger portion of stress and try to move away from one another (\cite{Kowalczyk-PreventingBlowUp-JMAA05}), whereas sensitivities of the type $S(n)\simeq\frac{1}{(n+1)^\alpha}$ can be motivated by the fact that movement in densely packed areas is inhibited by the large amount of present cells (\cite{HP-volumefilling-CAMQ02}).

Accordingly, extensive research has been committed to the study of \eqref{KS} with different varieties of $D(n)$ and $S(n)$ and their respective necessary conditions for global (and bounded) solutions to exist. An overview of different variations of the model and on the vast mathematical background can be found in the surveys \cite{HP09,BBWT15} and references therein. As one consequence of a long list of studies, from which we will only name a few and refer to the references in \cite{TaoWin-quasilinear_JDE12} for a more exhaustive overview, the critical exponent in the growth ratio of $\frac{S(n)}{D(n)}$ has been identified to be $\frac{2}{N}$. In fact, under the assumption of uniform parabolicity it was shown in \cite{TaoWin-quasilinear_JDE12} for the corresponding Neumann boundary value problem in a smooth domain $\Omega\subset\R^N$, that for any suitably regular initial data the classical solutions of \eqref{KS} are global and bounded whenever
\begin{align*}
\frac{S(n)}{D(n)}\leq C(n+1)^\beta\quad\text{for all }n\geq0\text{ with some }C>0\text{ and }\beta<\frac{2}{N},
\end{align*}
and in \cite{Win-volume-filling-MMAS10} the existence of smooth solutions which blow-up in either finite or infinite time has been proven in the case of
\begin{align*}
\frac{S(n)}{D(n)}\geq C n^\gamma\quad\text{for all }n>1\text{ with some }C>0\text{ and }\gamma>\frac{2}{N}.
\end{align*}
In particular, for the explicit case involving nondegenerate diffusion of porous medium type, i.e. $D(n)\equiv m(n+1)^{m-1}$, and $S(n)\equiv (1+n)^{1-\alpha}$, this leads to the condition $\alpha+m>\frac{2N-2}{N}$ for global solutions to exist (see also \cite{CieslakStinner-JDE15} for a result on finite time blow-up). 
\\[0.1cm]
\noindent{\textbf{Fluid interaction.}
In comparison, results for models incorporating fluid-interaction are less complete. The substantial effect fluid-interaction can have on the migration process is indicated by studies on broadcast spawning (e.g. \cite{coll1994chemical,miller1985demonstration}) or by the experiments undertaken in \cite{tuval2005bacterial}, where spontaneous turbulence emergence was observed with aerobic bacteria suspended in sessile drops of water. Since the bacteria consume the chemical instead of producing it, the authors of \cite{tuval2005bacterial} proposed the model
\begin{align}\label{CTcons}
\left\{
\begin{array}{r@{\,}c@{\,}c@{\ }l@{\quad}l}
n_{t}&+&u\cdot\!\nabla n&=\Delta n-\nabla\!\cdot(nS(c)\nabla c),\\
c_{t}&+&u\cdot\!\nabla c&=\Delta c-nc,\\
u_{t}&+&\kappa(u\cdot\nabla)u&=\Delta u+\nabla P+n\nabla\phi,\\
&&\nabla\cdot u&=0,
\end{array}\right.
\end{align}
for the unknown quantities $n,c,u,P$ denoting bacterial density, chemical concentration, fluid velocity and associated pressure, respectively, and $\phi$ is a prescribed gravitational potential. Apart from the biological motivation featuring aerobic bacteria, the consumption setting also has the minor advantage, that in contrast to its actual Keller--Segel-fluid counterpart (see \eqref{CTprod} below) one can immediately obtain uniform bounds on $c$ from the second equation, which led to it being studied more heavily than the framework with signal production by the cells. Let us briefly summarize some of the results available for \eqref{CTcons} in a three-dimensional bounded domain $\Omega$ with smooth boundary.

In the framework involving Stokes fluid (i.e. $\kappa=0$) and linear diffusion and where in fact $S=S(x,n,c)$ may be a tensor-valued function, accounting for modeling approaches where the direction of movement does not have to be parallel to the concentration gradient (\cite{OH-difflimit-SIAM02}), the existence of global classical solutions was shown for $\alpha>\frac{1}{6}$ (\cite{CaoWang-GlobClass-DCDSB15}), whereas in the full Navier--Stokes setting and consumption of the form $-nf(c)$, with some nonnegative $\CSp{1}{[0,\infty)}$ function $f$ satisfying $f(0)=0$ instead of just $-nc$ in the second equation, it was shown in \cite{win_globweak3d-AHPN16} that for arbitrary sufficiently smooth initial data there exist global weak solutions, whenever the scalar sensitivity function $S\in\CSp{2}{[0,\infty)}$ satisfies the hypotheses $(\frac{f}{S})'>0$, $(\frac{f}{S})''\leq0$ and $(S\cdot f)''\geq0$ on $[0,\infty)$. In the corresponding Cauchy-Problem similar results can be proven, as witnessed by the works \cite{DuanXiang-NoteOn-IMRN14,HeZhang-GlobEx-NARWA17}. The result for bounded domains has later also been extended to nonlinear diffusion of porous medium type $D(n)=mn^{m-1}$ for $m>\frac{2}{3}$ under the same conditions on $f$ and $S$ (\cite{ZhangLi-GlobWeakSol-JDE15}). In \cite{Vorotnikov-WeakSol-CMS14} a related system involving an additional source function $g(n)$ is studied and for $m\geq1$ and $g$ satisfying some growth conditions global weak solutions are obtained.

In the Stokes setting with $D(n)\equiv mn^{m-1}$ and $S(c)\equiv c$, the authors of \cite{TaoWin-LocBddGlobSol-AnnInstHP13} proved the existence of global weak solutions, which are locally bounded for any $m>\frac{8}{7}$. In \cite{Win-GlobExStabStrongDiff-17} one of the authors extended the global existence of weak solutions to values $m>\frac{9}{8}$ and also discussed convergence of these solutions to the spatial homogeneous steady state $(\frac{1}{|\Omega|}\intomega n_0,0,0)$.
One result concerning interplay between porous medium type diffusion and tensor-valued sensitivity satisfying $|S(x,n,c)|\leq\frac{S_0(c)}{(1+n)^\alpha}$ with some nondecreasing $S_0:[0,\infty)\to[0,\infty)$, is given in \cite{WangLi-ZAMP17}, which establishes the global existence of weak solutions for $m+\alpha>\frac{7}{6}$ and also verifies their convergence towards the steady state mentioned above. To conclude this nonexhaustive list, we mention the recent work \cite{ZhengLiZouZhangYan-17}, where global weak solutions in the Navier-Stokes setting with nonlinear diffusion $D\in\CSp{\gamma}{[0,\infty)}$ satisfying $C_1 n^{m-1}\geq D(n)\geq C_0n^{m-1}$ with $C_1,C_0>0$ and $m>\frac{10}{9}$ and a tensor-valued sensitivity satisfying $|S(x,n,c)|\leq S_0(c)$ with some nondecreasing $S_0:[0,\infty)\to\R$ were obtained. 

Concerning the framework where the chemical is produced by the cells instead of consumed, as in the actual Keller--Segel model, that is
\begin{align}\label{CTprod}
\left\{
\begin{array}{r@{\,}c@{\,}c@{\ }l@{\quad}l}
n_{t}&+&u\cdot\!\nabla n&=\nabla\!\cdot(D(n)\nabla n)-\nabla\!\cdot(nS(x,n,c)\nabla c),\\
c_{t}&+&u\cdot\!\nabla c&=\Delta c-c+n,\\
u_{t}&+&\kappa(u\cdot\nabla)u&=\Delta u+\nabla P+n\nabla\phi,\\
&&\nabla\cdot u&=0,
\end{array}\right.
\end{align}
where $S$ may again be a tensor-valued function, only a handful of results are available. On one hand, in a three-dimensional setup involving linear diffusion (i.e. $D(n)\equiv1$) and tensor-valued sensitivity $S(x,n,c)$ satisfying $|S(x,n,c)|\leq S_0(1+n)^{-\alpha}$ global weak solutions have been shown to exists in \cite{LiuWang-GlobWeak-JDE17} for $\alpha>\frac{3}{7}$ and global very weak solutions were obtained for $\alpha>\frac{1}{3}$ in \cite{Wang17-globweak-ksns}, which in light of the known results for the fluid-free system mentioned above is an optimal restriction on $\alpha$.
On the other hand, in a setting with diffusion of porous medium type (i.e. $D(n)\equiv mn^{m-1}$) and sensitivity $S(x,n,c)\equiv 1$ the only result we are aware of accounts for global weak solutions whenever $m>2$ (\cite{ZhengJ-globweak-JDE17}), which most probably is not optimal in the sense of $m+\alpha>\frac{2N-2}{N}$.

Furthermore, global existence for exponents smaller than the critical one $\alpha+m=\frac{2N-2}{N}$ can be obtained by including a logistic growth term of the form $+rn-\mu n^2$ in the first equation, as e.g. illustrated by the studies in \cite{LiuWang-GlobWeakLogSource-JDE16}, where $m+\alpha>\frac{6}{5}$ is sufficient to obtain global weak solutions in the three-dimensional Stokes setting (i.e. $\kappa=0$). 
\\[0.1cm]
\noindent{\textbf{Main results.} 
In a setting combining porous medium type diffusion and Navier--Stokes-fluid-interaction we attempt to attain optimal conditions on the diffusion exponent leading to global existence and therefore consider the prototypical system
\begin{align}\label{CTnonlindiff}
\left\{
\begin{array}{r@{\,}c@{\,}c@{\ }l@{\quad}l@{\quad}l@{\,}c}
n_{t}&+&u\cdot\!\nabla n&=\Delta n^m-\nabla\!\cdot(n\nabla c),\ &x\in\Omega,& t>0,\\
c_{t}&+&u\cdot\!\nabla c&=\Delta c-c+n,\ &x\in\Omega,& t>0,\\
u_{t}&+&(u\cdot\nabla)u&=\Delta u+\nabla P+n\nabla\phi,\ &x\in\Omega,& t>0,\\
&&\nabla\cdot u&=0,\ &x\in\Omega,& t>0,
\end{array}\right.
\end{align}
with boundary conditions
\begin{align}\label{BC}
\big(\nabla n^m(x,t)-n(x,t)\nabla c(x,t)\big)\cdot\nu=\nabla c(x,t)\cdot\nu=0\quad\text{and}\quad u(x,t)=0\qquad\text{for }x\in\romega\text{ and }t>0,
\end{align}
and initial conditions
\begin{align}\label{IC}
n(x,0)=n_0(x),\quad c(x,0)=c_0(x),\quad u(x,0)=u_0(x),\quad x\in\Omega,
\end{align}
where $\Omega\subset\R^3$ is a bounded domain with smooth boundary, $m>1$ and
\begin{align}\label{phidef}
\phi\in\W[2,\infty].
\end{align}
Assuming the initial data to satisfy
\begin{align}\label{IR}
\left\{\begin{array}{r@{\,}l}
n_0&\in\CSp{\gamma}{\bomega}\quad\text{for some }\gamma>0\quad\text{with } n_0\geq0\text{ in }\Omega\text{ and }n_0\not\equiv0,\\
c_0&\in\W[1,\infty]\quad\text{ with }c_0\geq0\text{ in }\bomega,\text{ and }c_0\not\equiv0,\\
u_0&\in\WSp{2,2}{\Omega;\R^3}\cap W_0^{1,2}\left(\Omega;\R^3\right)\quad\text{such that } \nabla\cdot u_0=0,
\end{array}\right.
\end{align}
we can state our main results as follows.

\begin{theorem}\label{theo:1}
Let $\Omega\subset\R^3$ be a bounded domain with smooth boundary. Suppose that $m>\frac{5}{3}$ and that $n_0,c_0$ and $u_0$ comply with \eqref{IR}. Then \eqref{CTnonlindiff}--\,\eqref{IC} admits at least one global weak solution in the sense of Definition \ref{def:weak-solution} below.
\end{theorem}

This extends the previous result of (\cite{ZhengJ-globweak-JDE17}), where the condition $m>2$ was obtained. On the other hand, for values of $m\in(\tfrac43,\tfrac53]$ we will only obtain compactness properties which seem too mild to pass to the limit in the weak formulation of our approximating system. Nevertheless, a very weak solution concept, which has been utilized in similar works before and is specified in Definition \ref{def:very_weak_sol} below, can still be administered with this weaker information, since in particular  $n\nabla n^{m-1}$ and $n\nabla c$ are not required to be integrable therein.

\begin{theorem}\label{theo:2}
Let $\Omega\subset\R^3$ be a bounded domain with smooth boundary. Suppose that $m>\frac{4}{3}$ and that $n_0,c_0$ and $u_0$ comply with \eqref{IR}. Then \eqref{CTnonlindiff}--\,\eqref{IC} admits at least one global very weak solution $(n,c,u)$ in the sense of Definition \ref{def:very_weak_sol} below. In particular, this global very weak solution satisfies
\begin{align*}
n&\in \LSploc{2m-\frac{4}{3}}{\bomega\times[0,\infty)},\quad c\in \LSploc{2}{[0,\infty);\W[1,2]},\quad u\in \LSploc{2}{[0,\infty);W_{0,\sigma}^{1,2}(\Omega;\R^3)},
\end{align*}
and
\begin{align*}
\intomega n(\cdot,t)=\intomega n_0\quad\text{for a.e. }t>0.
\end{align*}
\end{theorem}

Recalling the condition $m+\alpha>\frac{2N-2}{N}$ for global existence in the fluid-free setting, as implied by the previously mentioned studies \cite{TaoWin-quasilinear_JDE12,Win-volume-filling-MMAS10,CieslakStinner-JDE15}, this result appears to be optimal with respect to $m$.

\noindent{\textbf{Plan of the paper.} 
As our interest is mainly with small values of $m$, the main objective of our analysis will be to obtain a priori estimates, which capture optimal conditions on $m$. The fluid-coupling however destroys the well-known energy-structure of the classical Keller--Segel model and working with the standard energy functional cannot be expected to be of any help in deriving optimal a priori estimates in the setting of \eqref{CTnonlindiff}. While this difficulty can be circumvented in presence of a signal consuming process as in \eqref{CTcons}, by utilizing a suitable testing procedure to (more or less) cancel out the bad parts of the cross-diffusive term and obtain a quasi energy estimate, adapting such a testing procedure to the signal production present in \eqref{CTnonlindiff} seems rather hopeless when asking for optimal conditions on $m$, meaning that most sensible testing procedures which would improve the regularity information for $n$ beyond the obvious $L^1$-estimates are out of reach for small values of $m$. To counteract the missing energy estimate we will therefore investigate a functional of the form $\intomega n^{m-1}+\intomega c^2$, which for $m<2$ is obviously of sublinear growth in $n$ (which is rather uncommon functional to investigate) and seemingly does not improve our knowledge on the regularity of $n$, however, as byproduct of the underlying testing procedure we obtain a spatio-temporal estimate on the gradient term $\nabla n^{m-1}$ (cf. Lemma \ref{lem:bounds}) which, by means of standard Gagliardo--Nirenberg estimates, can be refined into a bound on $n\in\LSploc{2m-\frac43}{\bomega\times[0,\infty)}$ (cf. Lemma \ref{lem:st-bound-nepsi}), while still only prescribing the condition $m>\frac43$, and thereby (even for $m\in(\tfrac43,2)$) slightly improves the regularity information beyond the mere $L^1$-estimate. The very weak solution concept, specified in Definition \ref{def:very_weak_sol}, is mild enough to work with this minimum of regularity properties, while still being consistent with the concept of classical solutions. Assuming larger values of $m$, i.e. $m>\frac{5}{3}$, this analysis even provides sufficient regularity estimates for $n$ and $c$ to conclude that the integrals appearing in the weak solution concept remain finite. To specify, whenever $m>\frac{5}{3}$ we have $2m-\frac{4}{3}>2$ and hence $n\in\LSploc{2}{\bomega\times[0,\infty)}$, which upon combination with our other bounds suffices to obtain that $n\nabla n^{m-1}$, $n\nabla c$ and $nu$ belong to $\LSploc{1}{\bomega\times[0,\infty);\R^3}$ (cf. Lemma \ref{lem:weak-sol-5over3}).

Our setup is as follows: Starting with a brief introduction of the solution concepts we are going to consider (Sec. \ref{sec2:soldef}), we turn to a family of approximating systems allowing for global solutions (Sec. \ref{sec3:approxprob}) and discuss the previously mentioned a priori estimates (Sec. \ref{sec4:reg-est}) netting the cornerstone of our limit procedure (Sec. \ref{sec5:convergence}). Finally, depending on the size of $m$, we address the solution properties of the obtained limit functions (Sec. \ref{sec6:sol-prop}).

\setcounter{equation}{0} 
\section{Concepts of weak and very weak solvability}\label{sec2:soldef}
Since we  have two different concepts of solvability in the theorems above, in order to better differentiate between weak and very weak solutions, let us first specify what the very weak solutions we will obtain in Theorem \ref{theo:2} are supposed to satisfy. The concept draws on ideas originating from \cite{win15_chemorot} and \cite{Wang17-globweak-ksns}, which in our context has to be adapted to the nonlinear diffusion present in \eqref{CTnonlindiff}. The main difference to the standard notion of weak solvability lies in the fact that the first component is only expected to satisfy a kind of global supersolution property in the following sense.
\begin{definition}\label{def:weak_super_sol}
Let $\Phi\in\CSp{2}{[0,\infty)}$ be a nonnegative function satisfying $\Phi'>0$ on $(0,\infty)$. 
Assume that $n_0\in\Lo[\infty]$ is nonnegative and that $\Phi(n_0)\in\Lo[1]$. Moreover, let $c\in\LSploc{2}{[0,\infty);\W[1,2]}$ and $u\in\LSploc{1}{[0,\infty); W_0^{1,1}\left(\Omega;\R^3\right)}$ with $\dive u\equiv 0$ in $\mathcal{D}'\big(\Omega\times(0,\infty)\big)$. 
The nonnegative measurable function $n:\Omega\times(0,\infty)\to\R$ satisfying $n^{m-1}\in\LSploc{1}{[0,\infty);\W[1,1]}$ will be named a global weak $\Phi$--supersolution of the initial-boundary value problem
\begin{align}\label{eq:system_very_weak_supersol}
\left\{
\begin{array}{r@{\,}c@{\,}r@{\ }l@{\quad}l@{\quad}l@{\,}c}
n_{t}&+&u\cdot\!\nabla n&=\Delta n^m-\nabla\!\cdot(n\nabla c),\ &x\in\Omega,& t>0,\\
&&\frac{\partial n}{\partial \nu}&=0,\ &x\in\romega,& t>0,\\
&&n(x,0)&=n_0(x),\ &x\in\Omega,&
\end{array}\right.
\end{align}
if 
\begin{align}\label{eq:very_weak_supersol_regularity}
\Phi(n),\text{ and }\ \Phi''(n)n^{m-1}|\nabla n|^2&\text{ belong to }\LSploc{1}{\bomega\times[0,\infty)},\nonumber\\
\Phi'(n)n^{m-1}\nabla n,\text{ and }\ \Phi(n)u&\text{ belong to }\LSploc{1}{\bomega\times[0,\infty);\R^3},\\
\Phi'(n)n\text{ belongs to }\LSploc{2}{\bomega\times[0,\infty)},\text{ and }\ \Phi''(n)&n\nabla n\text{ belongs to }\LSploc{2}{\bomega\times[0,\infty);\R^3},\nonumber
\end{align}
and if for each nonnegative $\varphi\in C_0^\infty\left(\bomega\times[0,\infty)\right)$ with $\frac{\partial\varphi}{\partial\nu}=0$ on $\romega\times(0,\infty)$, the inequality
\begin{align}\label{eq:very_weak_supersol}
-\intinfomega&\Phi(n)\varphi_t-\intomega\Phi(n_0)\varphi(\cdot,0)\nonumber\\
&\geq - m\intinfomega\Phi''(n)n^{m-1}|\nabla n|^2\varphi-m\intinfomega\Phi'(n)n^{m-1}(\nabla n\cdot\nabla \varphi)\\&\qquad\qquad\,+\intinfomega\Phi''(n)n(\nabla n\cdot\nabla c)\varphi+\intinfomega\Phi'(n)n(\nabla c\cdot\nabla\varphi)+\intinfomega\Phi(n)(u\cdot\nabla\varphi)\nonumber
\end{align}
is satisfied.
\end{definition}

Later, for $m\in(\tfrac43,2)$, we will choose $\Phi(s)\equiv (s+1)^{m-1}$, for which $\Phi'(s)=(m-1)(s+1)^{m-2}$ and $\Phi''(s)=-(m-1)(2-m)(s+1)^{m-3}$. However, due to $m\in(\tfrac43,2)$ these quantities can be controlled from above by $(m-1)s^{m-2}$ and $-(m-1)(2-m)s^{m-3}$, respectively, for all $s\geq0$. Therefore, one of our main objectives will be to obtain bounds which will let us conclude that $n^{m-1}\in\LSploc{2}{[0,\infty);\W[1,2]}$, which (assuming $c$ and $u$ to be suitably regular) will suffice to treat all the integrals appearing in the supersolution property \eqref{eq:very_weak_supersol} (see also Corollary \ref{cor:1}, as well as \eqref{eq:sol-prop-n-eq} below). As for the other subproblems of \eqref{CTnonlindiff} we will require the properties for standard weak solvability to assemble the notion of global very weak solutions.

\begin{definition}\label{def:very_weak_sol}
A triple $(n,c,u)$ of functions
\begin{align*}
\begin{array}{r@{\ }c@{\ }l}
n&\in&\LSploc{1}{\bomega\times[0,\infty)},\\
c&\in&\LSploc{2}{[0,\infty);\W[1,2]},\\
u&\in& L_{loc}^1\big([0,\infty); W_{0}^{1,1}\!\left(\Omega;\R^3\right)\!\big),
\end{array}
\end{align*}
satisfying $n\geq0$ and $c\geq0$ in $\bomega\times[0,\infty)$, $cu\in\LSploc{1}{\bomega\times[0,\infty)}$, as well as $u\otimes u\in\LSploc{1}{\bomega\times[0,\infty);\R^{3\times3}}$ will be called a global very weak solution of \eqref{CTnonlindiff}--\,\eqref{IC}, if 
\begin{align*}
\intomega n(\cdot,t)\leq\intomega n_0\quad\text{for a.e. }t>0,
\end{align*}
if $\nabla\cdot u=0$ in $\mathcal{D}'\big(\Omega\times(0,\infty)\big)$, if the equality
\begin{align}\label{eq:very_weak_sol1}
-\intinfomega c\varphi_t-\intomega c_0\varphi(\cdot,0)=-\intinfomega\nabla c\cdot\nabla \varphi-\intinfomega c\varphi+\intinfomega n\varphi+\intinfomega c(u\cdot\nabla\varphi)
\end{align}
holds for all $\varphi\in\LSp{\infty}{\Omega\times(0,\infty)}\cap\LSp{2}{(0,\infty);\W[1,2]}$ with $\varphi_t\in\LSp{2}{\Omega\times(0,\infty)}$, which are compactly supported in $\bomega\times[0,\infty)$, if
\begin{align}\label{eq:very_weak_sol2}
-\intinfomega u\cdot\psi_t-\intomega u_0\cdot\psi(\cdot,0)=-\intinfomega\nabla u\cdot\nabla\psi+\intinfomega(u\otimes u)\cdot\nabla\psi+\intinfomega n\nabla\phi\cdot\psi
\end{align}
is fulfilled for all $\psi\in C_0^\infty\left(\Omega\times[0,\infty);\R^3\right)$ with $\nabla\cdot\psi\equiv0$ in $\Omega\times(0,\infty)$, and if finally there exists some nonnegative $\Phi\in\CSp{2}{[0,\infty)}$ with $\Phi'>0$ on $(0,\infty)$ such that $n$ is a global weak $\Phi$--supersolution of \eqref{eq:system_very_weak_supersol} in the sense of Definition \ref{def:weak_super_sol}.
\end{definition}

In contrast we will also talk about global weak solutions of \eqref{CTnonlindiff} in the standard sense, by which we mean the following.

\begin{definition}\label{def:weak-solution}
A triple $(n,c,u)$ of functions
\begin{align*}
n&\in\LSploc{1}{\bomega\times[0,\infty)},\\
c&\in\LSploc{1}{[0,\infty);\W[1,2]},\\
u&\in L_{loc}^1\big([0,\infty); W_{0}^{1,1}\!\left(\Omega;\R^3\right)\!\big),
\end{align*}
satisfying $n\geq0$ and $c\geq0$ in $\bomega\times[0,\infty)$, $n^{m-1}\in \LSploc{1}{[0,\infty);\W[1,1]}$ and $cu\in\LSploc{1}{\bomega\times[0,\infty);\R^3}$, as well as $u\otimes u\in\LSploc{1}{\bomega\times[0,\infty);\R^{3\times3}}$ will be called a global weak solution of \eqref{CTnonlindiff}--\,\eqref{IC}, if $\nabla\cdot u=0$ in $\mathcal{D}'\big((\Omega\times(0,\infty)\big)$, if
\begin{align*}
n\nabla n^{m-1}\quad\text{and}\quad n\nabla c,\quad\text{as well as}\quad nu&\quad\text{belong to }\LSploc{1}{\bomega\times[0,\infty);\R^3},
\end{align*}
if equality \eqref{eq:very_weak_sol1} holds for all $\varphi\in\LSp{\infty}{\Omega\times(0,\infty)}\cap\LSp{2}{(0,\infty);\W[1,2]}$ with $\varphi_t\in\LSp{2}{\Omega\times(0,\infty)}$, which are compactly supported in $\bomega\times[0,\infty)$, if  \eqref{eq:very_weak_sol2} is fulfilled for all $\psi\in C_0^\infty\left(\Omega\times[0,\infty);\R^3\right)$ with $\nabla\cdot\psi\equiv0$ in $\Omega\times(0,\infty)$, and if finally for each $\varphi\in C_0^\infty\big(\bomega\times[0,\infty)\big)$ with $\frac{\partial\varphi}{\partial\nu}=0$ on $\romega\times(0,\infty)$, the equality
\begin{align}\label{eq:weak-sol-n}
-\!\intinfomega n\varphi_t-\!\!\intomega n_0\,\varphi(\cdot,0)=-\frac{m}{m-1}\intinfomega n \big(\nabla n^{m-1}\!\cdot\!\nabla\varphi\big)+\!\!\intinfomega n(\nabla c\!\cdot\!\nabla\varphi)+\!\!\intinfomega n(u\cdot\!\nabla\varphi)
\end{align}
is satisfied.
\end{definition}

\begin{remark} i) If \eqref{eq:very_weak_supersol} is satisfied for $\Phi(s)\equiv s$ with equality, then $(n,c,u)$ is a global weak solution of \eqref{CTnonlindiff} in the sense of Definition \ref{def:weak-solution}, which shows that every global weak solution is also a global very weak solution.

ii) If the global very weak solution $(n,c,u)$ satisfies the regularity properties $n,c\in\CSp{0}{\bomega\times[0,\infty)}\cap\CSp{2,1}{\bomega\times(0,\infty)}$ and $u\in\CSp{0}{\bomega\times[0,\infty);\R^3}\cap\CSp{2,1}{\bomega\times(0,\infty);\R^3}$, it can be checked that the solution is also a global classical solution, i.e. one can find $P\in\CSp{1,0}{\bomega\times(0,\infty)}$ such that $(n,c,u,P)$ solves \eqref{CTnonlindiff} in the classical sense. See \cite[Lemma 2.1]{win15_chemorot} for the arguments involved. 
\end{remark}
\setcounter{equation}{0} 
\section{Global solutions to a family of approximating problems}\label{sec3:approxprob}
Working directly with the degenerate diffusion, the (possibly) unbounded chemotactic sensitivity, and the convection term present in the Navier-Stokes equation, poses quite some difficulties. Accordingly, we will fall back to a family of approximating problems regularized in a fashion which allows us to obtain global solutions in a straightforward manner. In fact, for $\epsi\in(0,1)$ we will consider the problems
\begin{align}\label{approxprob}
\left\{
\begin{array}{r@{\,}c@{\, }c@{\,}l@{\quad}l@{\quad}l@{\,}c}
n_{\epsi t}&+&u_\epsi\cdot\!\nabla n_\epsi&=\nabla\cdot\big(m(n_\epsi+\epsi)^{m-1}\nabla n_\epsi-\frac{n_\epsi}{(1+\epsi n_\epsi)^3}\nabla c_\epsi\big),\ &x\in\Omega,& t>0,\\
c_{\epsi t}&+&u_\epsi \cdot\!\nabla c_\epsi&=\Delta c_\epsi-c_\epsi+n_\epsi,\ &x\in\Omega,& t>0,\\
u_{\epsi t}&+&(Y_\epsi u_\epsi\cdot\nabla)u_\epsi&=\Delta u_\epsi+\nabla P_\epsi+n_\epsi\nabla\phi,\ &x\in\Omega,& t>0,\\
&&\quad\nabla\cdot u_\epsi&=0,\ &x\in\Omega,& t>0,\\
&&\quad\ \partial_\nu n_\epsi=\partial_\nu c_\epsi&=0,\qquad\qquad\quad u_\epsi=0,\ &x\in\romega,& t>0,\\
&&\qquad n_\epsi(x,0)&=n_0(x),\quad c_\epsi(x,0)=c_0(x),\quad u_\epsi(x,0)=u_0(x),\ &x\in\Omega,&
\end{array}\right.
\end{align}
where $Y_\epsi$ denotes the Yosida approximation of the Stokes operator given by
\begin{align*}
Y_\epsi\varphi:=(1+\epsi A)^{-1}\varphi\quad\text{for }\epsi\in(0,1)\text{ and }\varphi\in L_\sigma^2(\Omega).
\end{align*}
\subsection{Local existence of approximating solutions and basic properties}\label{sec31:locex}
Let us start by ensuring time-local existence of classical solutions to \eqref{approxprob}, which, including a suitable extensibility criterion, can be attained by employing well-known fixed point arguments. Denoting by $A:=\mathcal{P}\Delta$ the Stokes operator with Helmholtz projection $\mathcal{P}$ from $\Lo[2]$ to the solenoidal subspace $L_\sigma^2(\Omega;\R^3):=\{\varphi\in\LSp{2}{\Omega;\R^3}\,\vert\,\nabla\cdot\varphi=0\}$ we obtain the following.
\begin{lemma}\label{lem:loc_ex}
Let $\Omega\subset\R^3$ be a bounded domain with smooth boundary, $\phi\in\W[2,\infty]$, $\vartheta>3$ and $m\geq1$.  Suppose that $n_0,c_0$ and $u_0$ comply with \eqref{IR}. Then for any $\epsi\in(0,1)$, there exists $\Tme\in(0,\infty]$ and a uniquely determined triple $(n_\epsi,c_\epsi,u_\epsi)$ of functions satisfying
\begin{align*}
n_\epsi&\in \CSp{0}{\bomega\times[0,\Tme)}\cap \CSp{2,1}{\bomega\times(0,\Tme)}, \\
c_\epsi&\in \CSp{0}{\bomega\times[0,\Tme)}\cap \CSp{2,1}{\bomega\times(0,\Tme)}\cap\CSp{0}{[0,\Tme);\W[1,\vartheta]},\\
u_\epsi&\in \CSp{0}{\bomega\times[0,\Tme);\R^3}\cap \CSp{2,1}{\bomega\times(0,\Tme);\R^3},
\end{align*}
which, together with some $P_\epsi\in \CSp{1,0}{\bomega\times(0,\Tme)}$, solve \eqref{approxprob} in the classical sense and fulfill $n_\epsi\geq0$ and $c_\epsi\geq0$ in $\bomega\times[0,\Tme)$, as well as
\begin{align}\label{eq:loc-ex-alt}
\text{either}\quad\Tme=\infty\quad\text{or}\quad
\limsup_{t\nearrow\Tme}\big(\|n_\epsi(\cdot,t)\|_{\Lo[\infty]}&+\|c_\epsi(\cdot,t)\|_{\W[1,\vartheta]}+\|A^\beta u_\epsi(\cdot,t)\|_{\Lo[2]}\big)=\infty\nonumber\\&\text{for all }\vartheta>3\text{ and }\beta\in\big(\tfrac{3}{4},1\big).
\end{align}
\end{lemma}

\begin{bew}
Adapting well-established fixed point arguments as e.g. employed in \cite[Lemma 2.1]{tao_winkler_chemohapto11-siam11}, \cite[Lemma 2.2]{Lan17-LocBddGlobSolNonlinDiff-JDE} and \cite[Lemma 2.1]{win_fluid_final} for related frameworks, one can readily verify the existence of a local-in-time classical solution which satisfies \eqref{eq:loc-ex-alt}. The nonnegativity of the first two components is an immediate consequence of the maximum principle (\cite[Thm. 7.1.9]{evans}).
\end{bew}

In straightforward fashion one can check the boundedness of the $\Lo[1]$--norms, which is common in most chemotaxis settings.

\begin{lemma}\label{lem:l1-bounds}
Suppose that $m\geq 1$ and that $n_0,c_0$ and $u_0$ satisfy \eqref{IR}. Then for any $\epsi\in(0,1)$ the classical solution $(n_\epsi,c_\epsi,u_\epsi)$ of \eqref{approxprob} fulfills
\begin{align*}
\intomega n_\epsi(\cdot,t)=\intomega n_0\quad\text{for all }t\in(0,\Tme)
\end{align*}
and
\begin{align*}
\intomega c_\epsi(\cdot,t)\leq \max\left\{\intomega n_0,\intomega c_0\right\}\quad\text{for all }t\in(0,\Tme).
\end{align*}
\end{lemma}
\begin{bew}
The first statement can be obtained in standard manner by simple integration of the respective equation in \eqref{approxprob}. The second assertion then follows from integration of the second equation and an ODE comparison argument (\cite[Thm. IX]{Walter-ODE-98}).
\end{bew}

\subsection{Global approximating solutions}\label{sec32:globex}
In this section we want to ensure that the time-local solutions obtained in Lemma \ref{lem:loc_ex} are in fact global solutions. For this, we will rely on a Moser-type iteration (see e.g. \cite[Lemma A.1]{TaoWin-quasilinear_JDE12} for a version fitting our framework). In order to start the iteration process though, we will need additional regularity estimates for $n_\epsi,c_\epsi,\nabla c_\epsi$ and $u_\epsi$, which may depend on $\epsi$. In a first step we will combine two suitable differential inequalities to improve on the known smoothness for $n_\epsi,c_\epsi$ and $u_\epsi$.
\begin{lemma}\label{lem:globex-testing}
Let $m\geq 1$ and assume that $(n_0,c_0,u_0)$ comply with \eqref{IR} and that $\beta\in(\frac{3}{4},1)$. Then for any $T\in(0,\Tme]$ with $T<\infty$ and any $\epsi\in(0,1)$ there exists a constant $C=C(T,\epsi)$ such that the classical solution $(n_\epsi,c_\epsi,u_\epsi)$ of \eqref{approxprob} satisfies 
\begin{align*}
\intomega n_\epsi^6(\cdot,t)+\intomega c_\epsi^6(\cdot,t)\leq C\quad\text{for all }t\in(0,T),
\end{align*}
as well as
\begin{align*}
\|A^\beta u_\epsi(\cdot,t)\|_{\Lo[2]}+\|u_\epsi(\cdot,t)\|_{\Lo[\infty]}\leq C\quad\text{for all }t\in(0,T).
\end{align*}
\end{lemma}

\begin{bew}
We let $\gamma:=\max\{m-1,6\}$. For fixed $\epsi\in(0,1)$, we make use of the first equation in \eqref{approxprob}, integration by parts and the fact that $\nabla \cdot u_\epsi=0$ in $\Omega\times(0,\Tme)$ to calculate
\begin{align*}
\frac{1}{\gamma}\frac{\intd}{\intd t}\intomega n_\epsi^\gamma&=\intomega n_\epsi^{\gamma-1}\nabla\cdot\Big(m(n_\epsi+\epsi)^{m-1}\nabla n_\epsi-\frac{n_\epsi}{(1+\epsi n_\epsi)^3}\nabla c_\epsi\Big)-\frac{1}{\gamma}\intomega\nabla\cdot(n_\epsi^\gamma u_\epsi)\\
&=-(\gamma-1)m\intomega (n_\epsi+\epsi)^{m-1}n_\epsi^{\gamma-2}|\nabla n_\epsi|^2+(\gamma-1)\intomega \frac{n_\epsi^{\gamma-1}}{(1+\epsi n_\epsi)^3}(\nabla n_\epsi\cdot\nabla c_\epsi)
\end{align*}
on $(0,T)$. Now, since $-(s+\epsi)^{m-1}\leq -s^{m-1}$ for all $s\geq0$, $\frac{s}{1+\epsi s}\leq \frac{1}{\epsi}$ for all $s\geq0$ and $\gamma-m+1\geq0$, as well as $m-\gamma+5\geq0$ by choice of $\gamma$, an application of Young's inequality shows
\begin{align}\label{eq:globex-testing-eq1}
\frac{1}{\gamma}\frac{\intd}{\intd t}\intomega n_\epsi^\gamma&\leq -\frac{(\gamma-1)m}{2}\intomega n_\epsi^{m+\gamma-3}|\nabla n_\epsi|^2+\frac{\gamma-1}{2m\epsi^{\gamma-m+1}}\intomega\frac{|\nabla c_\epsi|^2}{(1+\epsi n_\epsi)^{m-\gamma+5}}\nonumber\\&\leq-\frac{(\gamma-1)m}{2}\intomega n_\epsi^{m+\gamma-3}|\nabla n_\epsi|^2+\frac{\gamma-1}{2m\epsi^{\gamma-m+1}}\intomega|\nabla c_\epsi|^2\quad\text{on }(0,T).
\end{align}
In a similar fashion, we multiply the second equation of \eqref{approxprob} with $(c_\epsi+1)^{\gamma-1}$ and again using that $u_\epsi$ is divergence-free, we integrate by parts to obtain
\begin{align*}
\frac{1}{\gamma}\frac{\intd}{\intd t}\intomega(c_\epsi+1)^\gamma+(\gamma-1)\intomega (c_\epsi+1)^{\gamma-2}|\nabla c_\epsi|^2+\intomega c_\epsi(c_\epsi+1)^{\gamma-1}=\intomega n_\epsi(c_\epsi+1)^{\gamma-1}
\end{align*}
on $(0,T)$, from which we infer by positivity of $c_\epsi$ and an application of Young's inequality that
\begin{align}\label{eq:globex-testing-eq2}
\frac{1}{\gamma m\epsi^{\gamma-m+1}}\frac{\intd}{\intd t}\intomega (c_\epsi+1)^{\gamma}+\frac{\gamma-1}{m\epsi^{\gamma-m+1}}\intomega|\nabla c_\epsi|^2\leq \frac{1}{\gamma m\epsi^{\gamma-m+1}}\intomega n_\epsi^\gamma+\frac{\gamma-1}{\gamma m\epsi^{\gamma-m+1}}\intomega(c_\epsi+1)^\gamma
\end{align}
holds on $(0,T)$. Thus, combining \eqref{eq:globex-testing-eq1} and \eqref{eq:globex-testing-eq2} and integrating the resulting inequality implies the existence of $C_1:=C_1(T,\epsi)$ satisfying
\begin{align}\label{eq:globex-testing-eq3}
\intomega n_\epsi^\gamma(\cdot,t)+\intomega (c_\epsi(\cdot,t)+1)^\gamma\leq C_1\quad\text{for all }t\in(0,T),
\end{align}
and thereby proves the first part of the lemma in light of the fact that $\gamma\geq6$. For the second part we first note that due to the continuous embedding $D(A^\beta)\hookrightarrow\CSp{\theta}{\bomega}$ for any $\theta\in(0,2\beta-\tfrac{3}{2})$ (see \cite[Lemma III.2.4.3]{sohr} and \cite[Thm. 5.6.5]{evans}), we only have to find $C_2>0$ such that $\|A^\beta u_\epsi(\cdot,t)\|_{\Lo[2]}\leq C_2$ holds for $t\in(0,T)$. For this, we first test the third equation of \eqref{approxprob} by $u_\epsi$ to obtain 
\begin{align}\label{eq:globex-testing-eq4}
\frac{1}{2}\frac{\intd}{\intd t}\intomega |u_\epsi|^2+\intomega|\nabla u_\epsi|^2=\intomega n_\epsi u_\epsi\cdot\nabla\phi\quad\text{for all }t\in(0,T),	
\end{align}
where we used the facts that $\nabla\cdot u_\epsi\equiv0$ and $\nabla \cdot(1+\epsi A)^{-1}u_\epsi\equiv0$. In light of \eqref{phidef} and \eqref{eq:globex-testing-eq3} this readily implies $\|u_\epsi(\cdot,t)\|_{\Lo[2]}\leq C_3$ in $(0,T)$ for some $C_3>0$. Relying on properties of the Yosida approximation $Y_\epsi$, we can also immediately find $C_4>0$ (cf. \cite[p.462 (3.6)]{MiyakawaSohr-MathZ88}) such that $v_\epsi:=(1+\epsi A)^{-1}u_\epsi$ satisfies
\begin{align*}
\|v_\epsi(\cdot,t)\|_{\Lo[\infty]}=\|(1+\epsi A)^{-1}u_\epsi(\cdot,t)\|_{\Lo[\infty]}\leq C_4\|u_\epsi(\cdot,t)\|_{\Lo[2]}\leq C_5:=C_3C_4\quad\text{for all }t\in(0,T).
\end{align*}
Finally, we can refine these bounds into the desired estimate for $\|A^\beta u_\epsi(\cdot,t)\|_{\Lo[2]}$ by a two-step procedure (see e.g. \cite[Lemma 3.9]{win_globweak3d-AHPN16}) by first testing the equation $u_{\epsi t}+A u_\epsi=\mathcal{P}(-(v_\epsi\cdot\nabla)u_\epsi+n_\epsi\nabla\phi)$ by $A u_\epsi$ netting $C_6>0$ such that
\begin{align*}
\intomega |\nabla u_\epsi|^2=\intomega |A^{\frac{1}{2}}u_\epsi|^2\leq C_6\quad\text{for all }t\in(0,T),
\end{align*}
and $C_7>0$ satisfying
\begin{align*}
\|\mathcal{P}((v_\epsi\cdot\nabla)u_\epsi+n_\epsi\nabla\phi)\|_{\Lo[2]}\leq C_7\quad\text{for all }t\in(0,T).
\end{align*}
Secondly, we express $A^\beta u_\epsi$ by its variation-of-constants representation and make use of well-known smoothing properties of the Stokes semigroup (e.g. \cite[Lemma 3.1]{Win-ct_fluid_3d-CPDE15}) to obtain $C_8>0$ such that
\begin{align*}
\|A^\beta u_\epsi(\cdot,t)\|_{\Lo[2]}\leq C_8t^{-\beta}\|u_0\|_{\Lo[2]}+\frac{C_8 T^{1-\beta}}{1-\beta}\quad\text{for all }t\in(0,T),
\end{align*}
which completes the proof.
\end{bew}

The lemma above at hand, we can now obtain information on the gradient of $c_\epsi$, which will be the essential ingredient in order to satisfy the requirements of the Moser-type iteration, from which we will conclude that for each $\epsi\in(0,1)$ we have $\Tme=\infty$.

\begin{lemma}\label{lem:globex}
Let $m\geq 1$ and suppose that $n_0,c_0$ and $u_0$ satisfy \eqref{IR} and that $\beta\in(\frac{3}{4},1)$. Then for all $\epsi\in(0,1)$ the solution $(n_\epsi,c_\epsi,u_\epsi)$ of \eqref{approxprob} satisfies $\Tme=\infty$.
\end{lemma}

\begin{bew}
As a preliminary step we will require some regularity on $\nabla c_\epsi$. For this we fix $\epsi\in(0,1)$, assume that $\Tme<\infty$ and test the second equation of \eqref{approxprob} by $-\Delta c_\epsi$ and obtain, upon two applications of Young's inequality, that
\begin{align*}
\frac{1}{2}\frac{\intd}{\intd t}\intomega |\nabla c_\epsi|^2+\intomega|\Delta c_\epsi|^2+\intomega|\nabla c_\epsi|^2&=-\intomega n_\epsi\Delta c_\epsi+\intomega (u_\epsi\cdot\nabla c_\epsi)\Delta c_\epsi\\
&\leq \intomega n_\epsi^2+\frac{1}{2}\intomega|\Delta c_\epsi|^2+\|u_\epsi(\cdot,t)\|_{\Lo[\infty]}^2\intomega|\nabla c_\epsi(\cdot,t)|^2
\end{align*}
holds on $(0,\Tme)$. Recalling the bounds provided by Lemma \ref{lem:globex-testing}, this immediately implies
\begin{align*}
\intomega |\nabla c_\epsi(\cdot,t)|^2\leq C_1\quad\text{for all }t\in(0,\Tme)
\end{align*}
with some $C_1>0$. Next, we can combine the bounds provided by Lemma \ref{lem:globex-testing} with the new information on the spatial gradient of $c_\epsi$ and well-known smoothing properties of the Neumann heat semigroup (e.g. \cite[Lemma 1.3]{win10jde}) to find $C_2>0$ such that 
\begin{align*}
\|\nabla c_\epsi(\cdot,t)\|_{\Lo[\frac{11}{2}]}\leq C_2\quad\text{for all }t\in(0,\Tme),
\end{align*}
by simple expression of $\nabla c_\epsi$ in its corresponding variation-of-constants representation. In fact we now have $\frac{n_\epsi(\cdot,t)}{(1+\epsi n_\epsi(\cdot,t))^3}\nabla c_\epsi(\cdot,t)+n_\epsi(\cdot,t) u_\epsi(\cdot,t)\in\Lo[q]$ for all $t\in(0,\Tme)$, with some $q>5$ and hence we may employ a Moser type iteration (see \cite[Lemma A.1]{TaoWin-quasilinear_JDE12} for a version applicable to our system) to find $C_3>0$ such that $\|n_\epsi(\cdot,t)\|_{\Lo[\infty]}\leq C_3$ holds for all $t\in(0,\Tme)$.

Now, we see that combining the bound for $A^\beta u_\epsi(\cdot,t)$ in $\Lo[2]$, as contained in Lemma \ref{lem:globex-testing}, with the bounds prepared in the first part of this proof entails the existence of $C_4>0$ satisfying
\begin{align*}
\|n_\epsi(\cdot,t)\|_{\Lo[\infty]}+\|c_\epsi(\cdot,t)\|_{\W[1,5]}+\|A^\beta u_\epsi(\cdot,t)\|_{\Lo[2]}\leq C_4\quad\text{for all }t\in(0,\Tme),
\end{align*}
which, by our assumption of $\Tme<\infty$, clearly contradicts \eqref{eq:loc-ex-alt} and thereby proves $\Tme=\infty$.
\end{bew}

\setcounter{equation}{0} 
\section{Regularity estimates independent of \texorpdfstring{$\epsi$}{epsilon}}\label{sec4:reg-est}
Our main objective in this section will be to derive regularity information which is independent on $\epsi\in(0,1)$, while maintaining optimal conditions on $m$. Currently, the $\Lo[1]$--estimates present in Lemma \ref{lem:l1-bounds} are our only knowledge of this kind. Since we cannot rely on well-established testing procedures for the standard Keller-Segel system to improve the known information on $n_\epsi$, due to the fluid terms present in \eqref{CTnonlindiff}, we will investigate the functional $\intomega n_\epsi^{m-1}(\cdot,t)$, which for small values of $m>1$ is even of sublinear growth (cf. Lemma \ref{lem:testing-nepsi-m-1}). While at first (at least for $m<2$) this appears to not provide new information whatsoever, coupling this functional with $\intomega c_\epsi^2(\cdot,t)$ makes it possible to obtain a first information on the spatial gradient of $n_\epsi$ (cf. Lemma \ref{lem:bounds}), which in a second step can be refined to slightly more regularity information on $n_\epsi$ (cf. Lemma \ref{lem:st-bound-nepsi}). In the later parts of this section we then prepare all remaining bounds necessary for the limiting procedure undertanken in Section \ref{sec5:convergence}.

\subsection{Core estimates on the regularity of \texorpdfstring{$n_\epsi$ and $c_\epsi$}{the approximating solutions}}\label{sec41:regesti}
In preparation of some of our testing procedures we state the following elementary lemma.
\begin{lemma}\label{lem:testing-nepsi-m-1}
Let $m>1$ and assume that $n_0,c_0$ and $u_0$ comply with \eqref{IR}. Then for any $\epsi\in(0,1)$ and each $\varphi\in C^\infty\left(\bomega\times[0,\infty)\right)$ with $\frac{\partial\varphi}{\partial\nu}=0$ on $\romega\times(0,\infty)$ the classical solution $(n_\epsi,c_\epsi,u_\epsi)$ of \eqref{approxprob} satisfies
\begin{align}\label{eq:testing-nepsi}
\intomega \big((n_\epsi+\epsi)^{m-1}\big)_t\varphi=\frac{m(2-m)}{m-1}&\intomega |\nabla(n_\epsi+\epsi)^{m-1}|^2\varphi-m\intomega (n_\epsi+\epsi)^{m-1}\big(\nabla (n_\epsi+\epsi)^{m-1}\cdot\nabla \varphi\big)\nonumber\\
-(2&-m)\intomega\frac{n_\epsi(n_\epsi+\epsi)^{-1}}{(1+\epsi n_\epsi)^3}(\nabla (n_\epsi+\epsi)^{m-1}\cdot\nabla c_\epsi)\,\varphi\\\nonumber&\ +(m-1)\intomega \frac{n_\epsi(n_\epsi+\epsi)^{m-2}}{(1+\epsi n_\epsi)^3}(\nabla c_\epsi\cdot\nabla \varphi)+\intomega (n_\epsi+\epsi)^{m-1}(u_\epsi \cdot\nabla\varphi)
\end{align}
 on $(0,\infty)$.
\end{lemma}
\begin{bew}
In light of \eqref{approxprob} and the fact that $\nabla\cdot u_\epsi\equiv 0$ in $\Omega\times(0,\infty)$, we see that
\begin{align*}
\ &\intomega\big((n_\epsi+\epsi)^{m-1}\big)_t\varphi\\=\ &(m-1)\intomega (n_\epsi+\epsi)^{m-2}\varphi\,\nabla\cdot\Big(m(n_\epsi+\epsi)^{m-1}\nabla n_\epsi-\frac{n_\epsi}{(1+\epsi n_\epsi)^3}\nabla c_\epsi\Big)-\intomega \nabla\cdot\big((n_\epsi+\epsi)^{m-1}u_\epsi\big)\,\varphi
\end{align*}
holds for all $t>0$. Hence, the assertion follows from straightforward integration by parts and rewriting the resulting terms.
\end{bew}

In order to obtain any information on $\nabla n_\epsi^{m-1}$ whatsoever, we have to face the obstacle that $(2-m)$ is positive for small values of $m$. The key idea will be to employ Lemma \ref{lem:testing-nepsi-m-1} for a constant test function with negative sign, making it possible to transfer the term $\intomega|\nabla n_\epsi^{m-1}|^2$ to the left hand side of \eqref{eq:testing-nepsi}. Similar ideas have previously been used with success in e.g. \cite[Lemma 4.1]{Wang17-globweak-ksns}.

\begin{lemma}\label{lem:bounds}
Let $m>\frac{4}{3}$ and suppose that $n_0,c_0$ and $u_0$ fulfill \eqref{IR}. Then there exists some $C>0$ such that for all $\epsi\in(0,1)$ the global classical solution $(n_\epsi,c_\epsi,u_\epsi)$ of \eqref{approxprob} satisfies
\begin{align}\label{eq:lem-bounds-eq}
\intomega (n_\epsi+\epsi)^{m-1}(\cdot,t)+\intomega c_\epsi^2(\cdot,t) +\int_t^{t+1}\!\intomega\big|\nabla (n_\epsi+\epsi)^{m-1}\big|^2+\int_t^{t+1}\!\intomega|\nabla c_\epsi|^2\leq C
\end{align}
for all $t\geq0$.
\end{lemma}

\begin{bew} 
We will mainly concern ourselves with the case $m\in(\tfrac{4}{3},2)$ and give a few comments on necessary adjustments for the cases $m>2$ and $m=2$ at the end of the proof. For $m\in(\tfrac43,2)$ we employ Lemma \ref{lem:testing-nepsi-m-1} with $\varphi=-\frac{1}{m-1}$ to find that
\begin{align*}
-\frac{1}{m-1}\frac{\intd}{\intd t}\intomega (n_\epsi+\epsi)^{m-1}=-\frac{m(2-m)}{(m-1)^2}&\intomega \big|\nabla( n_\epsi+\epsi)^{m-1}\big|^2\\&\qquad+\frac{2-m}{m-1}\intomega\frac{n_\epsi(n_\epsi+\epsi)^{-1}}{(1+\epsi n_\epsi)^3}\big(\nabla(n_\epsi+\epsi)^{m-1}\cdot\nabla c_\epsi\big)\nonumber
\end{align*}
holds on $(0,\infty)$ for all $\epsi\in(0,1)$. Hence, making use of Young's inequality and the fact that for any $\epsi\in(0,1)$ we have $\frac{n_\epsi(n_\epsi+\epsi)^{-1}}{(1+\epsi n_\epsi)^3}\leq1$ in $\Omega\times(0,\infty)$, we obtain
\begin{align}\label{eq:bounds-proof-eq-n}
-\frac{1}{m-1}\frac{\intd}{\intd t}\intomega (n_\epsi+\epsi)^{m-1}(\cdot,t)&\leq-\frac{m(2-m)}{2(m-1)^2}\intomega \big|\nabla (n_\epsi+\epsi)^{m-1}(\cdot,t)\big|^2+\frac{2-m}{2m}\intomega|\nabla c_\epsi(\cdot,t)|^2
\end{align}
for all $\epsi\in(0,1)$ and all $t>0$. On the other hand, testing the second equation of \eqref{approxprob} by $c_\epsi$ we see that
\begin{align*}
\frac{1}{2}\frac{\intd}{\intd t}\intomega c_\epsi^2(\cdot,t)+\intomega|\nabla c_\epsi(\cdot,t)|^2+\intomega c_\epsi^2(\cdot,t)\leq \|c_\epsi(\cdot,t)\|_{\Lo[6]}\|n_\epsi(\cdot,t)\|_{\Lo[\nfrac{6}{5}]}
\end{align*}
is valid for all $\epsi\in(0,1)$ and all $t>0$ in light of Hölder's inequality and $u_\epsi$ being divergence-free. Making use of the embedding $\W[1,2]\hookrightarrow\Lo[6]$ and Young's inequality we thereby obtain $C_1>0$ such that
\begin{align}\label{eq:bounds-proof-eq-c}
\frac{\intd}{\intd t}\intomega c_\epsi^2(\cdot,t)+\intomega|\nabla c_\epsi(\cdot,t)|^2+\intomega c_\epsi^2(\cdot,t)\leq 2C_1^2\|n_\epsi(\cdot,t)\|_{\Lo[\nfrac{6}{5}]}^2\quad\text{for all }t>0\text{ and all }\epsi\in(0,1).
\end{align}
Combining \eqref{eq:bounds-proof-eq-n} with a multiple of \eqref{eq:bounds-proof-eq-c} we find $C_2:=\frac{2-m}{m}2C_1^2>0$ satisfying
\begin{align}\label{eq:bounds-proof-eq-comb}
\frac{\intd}{\intd t}\Big[-\frac{1}{m-1}&\intomega (n_\epsi+\epsi)^{m-1}(\cdot,t)+\frac{2-m}{m}\intomega c_\epsi^2(\cdot,t)\Big]+\frac{2-m}{m}\intomega c_\epsi^2(\cdot,t)\\&+\frac{m(2-m)}{2(m-1)^2}\intomega\big|\nabla (n_\epsi+\epsi)^{m-1}(\cdot,t)\big|^2+\frac{2-m}{2m}\intomega|\nabla c_\epsi(\cdot,t)|^2\leq C_2\|n_\epsi(\cdot,t)\|_{\Lo[\nfrac{6}{5}]}^2\nonumber
\end{align}
for all $t>0$ and all $\epsi\in(0,1)$. To further estimate the right hand side, we may employ the \GNI , Lemma \ref{lem:l1-bounds}, the nonnegativity of $n_\epsi$ and the fact that $\epsi<1$ to obtain $C_3>0$ such that
\begin{align*}
C_2\|n_\epsi\|_{\Lo[\nfrac{6}{5}]}^2\leq C_2\|n_\epsi+\epsi\|_{\Lo[\nfrac{6}{5}]}^2\leq C_3\|\nabla(n_\epsi+\epsi)^{m-1}\|_{\Lo[2]}^{\frac{2}{6m-7}}+C_3\quad\text{on }(0,\infty)\text{ for all }\epsi\in(0,1).
\end{align*}
Now, since $m>\frac{4}{3}$, clearly $\frac{2}{6m-7}<2$ and hence Young's inequality provides $C_4>0$ satisfying
\begin{align}\label{eq:bounds-proof-eq-n65}
C_2\|n_\epsi(\cdot,t)\|_{\Lo[\nfrac{6}{5}]}^2\leq \frac{m(2-m)}{4(m-1)^2}\intomega\big|\nabla(n_\epsi+\epsi)^{m-1}(\cdot,t)\big|^2+C_4\quad\text{for all }t>0\text{ and all }\epsi\in(0,1).
\end{align}
Consequently, letting 
\begin{align*}
y_\epsi(t):=-\frac{1}{m-1}\intomega (n_\epsi+\epsi)^{m-1}(\cdot,t)+\frac{2-m}{m}\intomega c_\epsi^2(\cdot,t),\quad t>0,
\end{align*}
and 
\begin{align*}
g_\epsi(t):=\frac{m(2-m)}{4(m-1)^2}\intomega\big|\nabla(n_\epsi+\epsi)^{m-1}(\cdot,t)\big|^2+\frac{2-m}{2m}\intomega|\nabla c_\epsi(\cdot,t)|^2,\quad t>0,
\end{align*}
we see by combination of \eqref{eq:bounds-proof-eq-comb} and \eqref{eq:bounds-proof-eq-n65} that in light of the fact that $y_\epsi(t)\leq \frac{2-m}{m}\intomega c_\epsi^2(\cdot,t)$ holds for all $t>0$ and all $\epsi\in(0,1)$, we have
\begin{align}\label{eq:bounds-eq-diffin}
y_\epsi'(t)+y_\epsi(t)+g_\epsi(t)\leq C_4\quad\text{for all }t>0\text{ and all }\epsi\in(0,1).
\end{align}
Since $g_\epsi\geq0$ for all $t>0$, an ODE comparison argument thereby implies that
\begin{align*}
y_\epsi(t)\leq C_5:=\max\Big\{-\frac{1}{m-1}\intomega (n_0+1)^{m-1}+\frac{2-m}{m}\intomega c_0^2,\ C_4\Big\}\quad\text{for all }t>0\text{ and all }\epsi\in(0,1),
\end{align*}
which does not imply the asserted bounds as of yet, since $y_\epsi(t)$ might in fact be negative. Nevertheless, since $m<2$ the claimed boundedness of $\intomega\!\;(n_\epsi+\epsi)^{m-1}$ is an immediate consequence of Lemma \ref{lem:l1-bounds} and hence there exists $C_6>0$ such that $\frac{1}{m-1}\intomega\;\!(n_\epsi+\epsi)^{m-1}\leq C_6$ for all $t>0$. Combining this with the estimate for $y_\epsi(t)$ we find that
\begin{align*}
\frac{2-m}{m}\intomega c_\epsi^2(\cdot,t)\leq \frac{1}{m-1}\intomega (n_\epsi+\epsi)^{m-1}(\cdot,t)+C_5\leq C_6+C_5\quad\text{holds for all }t>0\text{ and all }\epsi\in(0,1).
\end{align*}
As for the integral containing the derivatives in \eqref{eq:lem-bounds-eq}, we observe that \eqref{eq:bounds-eq-diffin} also shows that
\begin{align*}
\int_t^{t+1}g_\epsi(s)\intd s\leq y_\epsi(t)-y_\epsi(t+1)-\int_t^{t+1}y_\epsi(s)\intd s+C_4\quad\text{for all }t\geq0\text{ and all }\epsi\in(0,1),
\end{align*}
where by the definition of $y_\epsi$ and the positivity of $(2-m)\intomega c_\epsi^2(\cdot,t)$ for $t>0$, we may rely once more on Lemma \ref{lem:l1-bounds}  to estimate
\begin{align*}
-y_\epsi(t)\leq\frac{1}{m-1}\intomega (n_\epsi+\epsi)^{m-1}(\cdot,t)\leq C_6\quad\text{for all }t\geq 0\text{ and all }\epsi\in(0,1),
\end{align*}
so that in fact
\begin{align*}
\int_t^{t+1}g_\epsi(s)\intd s\leq C_5+2C_6+C_4\quad\text{for all }t\geq0\text{ and all }\epsi\in(0,1),
\end{align*}
proving the boundedness of the remaining integrals in \eqref{eq:lem-bounds-eq}.

To obtain the desired bound in the case of $m>2$, we repeat the steps above with $\varphi=\frac{1}{m-1}$ instead (see e.g. \cite[Lemma 2.3]{lankchapto15} for a version of the \GNI\ allowing for the $L^p$--spaces with $p<1$ required in this case), which upon combination with \eqref{eq:bounds-proof-eq-c} leads to a differential inequality of the kind featured in \eqref{eq:bounds-eq-diffin}, where this time the prefactor of $\intomega\!\;(n_\epsi+\epsi)^{m-1}$ in $y_\epsi(t)$ is positive, i.e.
\begin{align*}
\frac{\intd}{\intd t}\Big[\frac{1}{m-1}\intomega(n_\epsi+\epsi)^{m-1}(\cdot,t)&+\frac{m-2}{m}\intomega c^2_\epsi(\cdot,t)\Big]+\frac{m-2}{m}\intomega c_\epsi^2+\frac{m(m-2)}{4(m-1)^2}\intomega \big|\nabla(n_\epsi+\epsi)^{m-1}(\cdot,t)\big|^2\\&+\frac{m-2}{2m}\intomega|\nabla c_\epsi(\cdot,t)|^2\leq C_7\quad\text{for all }t>0\text{ and all }\epsi\in(0,1),
\end{align*}
with some $C_7>0$. Estimating the gradient term of $(n_\epsi+\epsi)^{m-1}$ from below by the \GNI\ in turn implies the asserted bound of $\intomega\;\!(n_\epsi+\epsi)^{m-1}$ and, due to the positivity of $y_\epsi(t)$ in this case, the conclusion of $\int_t^{t+1}g_\epsi(s)\intd s\leq C$ follows directly from the differential inequality and the bound for $y_\epsi(t)$. In the case of $m=2$ we estimate
\begin{align*}
\frac{\intd}{\intd t}\intomega (n_\epsi\ln n_\epsi)(\cdot,t)\leq -\intomega\big|\nabla (n_\epsi+\epsi)(\cdot,t)\big|^2+\frac{1}{4}\intomega|\nabla c_\epsi(\cdot,t)|^2\quad\text{for all }t>0,
\end{align*}
and combine with \eqref{eq:bounds-proof-eq-c} again to conclude the boundedness of the asserted integrals in a similar fashion as before, while making use of the fact that $s\ln s\geq -\frac{1}{e}$ for all $s>0$.
\end{bew}

With the latter spatio-temporal bound for $\nabla(n_\epsi+\epsi)^{m-1}$ at hand, we can now establish the following spatio-temporal bounds for $n_\epsi+\epsi$, which will play a key role in deriving uniform bounds for $u_\epsi$ and convergence properties for $n_\epsi$.

\begin{lemma}\label{lem:st-bound-nepsi}
Let $m>\frac{4}{3}$ and assume that $n_0,c_0$ and $u_0$ comply with \eqref{IR}. Then for all $p\in\big(1,6(m-1)\big)$ there exists $C>0$ such that for all $\epsi\in(0,1)$ the solution $(n_\epsi,c_\epsi,u_\epsi)$ of \eqref{approxprob} satisfies
\begin{align}\label{eq:st-bound-np}
\int_t^{t+1}\! \big\|n_\epsi(\cdot,s)+\epsi\big\|_{\Lo[p]}^{\frac{2p(m-\frac{7}{6})}{p-1}}\intd s\leq C\quad\text{for all }t>0.
\end{align}
In particular, there exists $C>0$ such that
\begin{align}\label{eq:st-bound-specialcases}
\int_t^{t+1}\!\big\|n_\epsi(\cdot,s)+\epsi\big\|_{\Lo[\frac{6}{5}]}^2\intd s\leq C\quad\text{and}\quad\int_t^{t+1}\big\|n_\epsi(\cdot,s)+\epsi\big\|_{\Lo[2m-\frac{4}{3}]}^{2m-\frac{4}{3}}\intd s\leq C
\end{align}
hold for each $\epsi\in(0,1)$ and all $t\geq0$.
\end{lemma}

\begin{bew} Inspired by the arguments of \cite[Lemma 4.2]{Wang17-globweak-ksns}, we employ the \GNI\ (see e.g. \cite[Lemma 2.3]{lankchapto15}) to obtain $C_1>0$ such that
\begin{align*}
\int_t^{t+1}\!\!\!\big\|n_\epsi&(\cdot,s)+\epsi\big\|_{\Lo[p]}^{\frac{2p(m-\frac{7}{6})}{p-1}}\intd s=\int_t^{t+1}\!\!\!\big\|(n_\epsi+\epsi)^{m-1}(\cdot,s)\big\|_{\Lo[\frac{p}{m-1}]}^{\frac{2p}{p-1}\cdot\frac{6m-7}{6(m-1)}}\intd s\\
&\leq C_1\!\int_t^{t+1}\!\!\!\big\|\nabla(n_\epsi+\epsi)^{m-1}(\cdot,s)\big\|_{\Lo[2]}^{\frac{2p}{p-1}\cdot\frac{6m-7}{6(m-1)}\cdot a}\big\|(n_\epsi+\epsi)^{m-1}(\cdot,s)\big\|_{\Lo[\frac{1}{m-1}]}^{\frac{2p}{p-1}\cdot\frac{6m-7}{6(m-1)}\cdot(1-a)}\intd s\\&\hspace*{7.3cm}+C_1\!\int_t^{t+1}\!\!\!\big\|(n_\epsi+\epsi)^{m-1}(\cdot,s)\big\|_{\Lo[\frac{1}{m-1}]}^{\frac{2p}{p-1}\cdot\frac{6m-7}{6(m-1)}}\intd s
\end{align*}
holds for all $t\geq0$ and all $\epsi\in(0,1)$, where
\begin{align*}
a=\frac{m-1-\frac{m-1}{p}}{m-1+\frac{1}{3}-\frac{1}{2}}=\frac{p-1}{p}\cdot\frac{6(m-1)}{6m-7}\in(0,1)
\end{align*}
due to $p\in(1,6(m-1))$ and $m>\frac{7}{6}$. In consideration of Lemma \ref{lem:l1-bounds} this entails the existence of $C_2>0$ satisfying
\begin{align*}
\int_t^{t+1}\!\big\|n_\epsi(\cdot,s)+\epsi\big\|_{\Lo[p]}^{\frac{2p(m-\frac{7}{6})}{p-1}}\intd s\leq C_2\int_t^{t+1}\!\intomega\big|\nabla (n_\epsi+\epsi)^{m-1}\big|^2+C_2\quad\text{for all }t\geq0\text{ and all }\epsi\in(0,1),
\end{align*}
which, due to $m>\frac{4}{3}$, immediately implies \eqref{eq:st-bound-np} in light of Lemma \ref{lem:bounds}. As for the special cases in \eqref{eq:st-bound-specialcases}, we only have to ensure that each of these $p$ satisfy $p\in(1,6(m-1))$ and that the given exponent is less than or equal to $\frac{2p(m-\frac{7}{6})}{p-1}$, since then, with the bound from the first step at hand, an application of Young's inequality directly implies the assertion. In both cases these conditions are fulfilled as an immediate consequence of the fact that $m>\frac{4}{3}$.
\end{bew}

Let us also briefly prepare some additional bounds, which will play an important role in the limit process for the explicit choice of $\Phi(s)=(s+1)^{m-1}$ with $m\in(\tfrac43,2)$.
\begin{corollary}\label{cor:1}
Let $m\in(\tfrac43,2)$ and suppose that $n_0,c_0$ and $u_0$ fulfill \eqref{IR}. Then there exists some $C_1>0$ such that for all $\epsi\in(0,1)$ the global classical solution $(n_\epsi,c_\epsi,u_\epsi)$ of \eqref{approxprob} satisfies
\begin{align}\label{eq:cor1}
\int_t^{t+1}\!\intomega\big|\nabla(n_\epsi+1)^{m-1}\big|^2+\int_t^{t+1}\!\intomega\big|(n_\epsi+1)^{\frac{m-3}{2}}(n_\epsi+\epsi)^{\frac{m-1}{2}}\nabla n_\epsi\big|^2\leq C_1,
\end{align}
for all $t\geq0$. Moreover, there exist $p>2$ and $C_2>0$ such that
\begin{align}\label{eq:cor2}
\int_t^{t+1}\big\|(n_\epsi+1)^{m-1}\big\|_{\Lo[p]}^p\leq C_2
\end{align}
holds for each $\epsi\in(0,1)$ and all $t\geq0$.
\end{corollary}
\begin{bew}
Since $1<m<2$, we can easily estimate
\begin{align*}
\int_t^{t+1}\!\intomega\big|\nabla(n_\epsi+1)^{m-1}\big|^2&\leq\int_t^{t+1}\!\intomega(m-1)^2(n_\epsi+\epsi)^{2(m-2)}\big|\nabla n_\epsi\big|^2=\int_t^{t+1}\!\intomega\big|\nabla(n_\epsi+\epsi)^{m-1}\big|^2
\end{align*}
for all $t\geq0$ and all $\epsi\in(0,1)$. Hence, the boundedness of the first term in \eqref{eq:cor1} is a direct consequence of Lemma \ref{lem:bounds}. Since $m>1$, the remaining bound in \eqref{eq:cor1} follows immediately from the one we just established.  For the second part we note that due to $m>\tfrac76$ the interval $(\max\{2,\tfrac{1}{m-1}\},\tfrac{6m-4}{3m-3})$ is not empty and hence, we can fix $p\in(\max\{2,\tfrac{1}{m-1}\},\tfrac{6m-4}{3m-3})$ and then employ the \GNI\ to find $C>0$ such that
\begin{align*}
\int_t^{t+1}\!\big\|(n_\epsi+1)^{m-1}\big\|_{\Lo[p]}^p\leq C\int_t^{t+1}\!\big\|\nabla (n_\epsi+1)^{m-1}\big\|_{\Lo[2]}^{pa}\big\|(n_\epsi&+1)^{m-1}\big\|_{\Lo[\frac1{m-1}]}^{p(1-a)}\\&+C\int_t^{t+1}\!\big\|(n_\epsi+1)^{m-1}\big\|_{\Lo[\frac1{m-1}]}^p
\end{align*}
for all $t\geq0$ and all $\epsi\in(0,1)$, where $a=\frac{6p(m-1)-6}{(6m-7)p}$. Our choice of $p$ implies $pa<2$ and therefore, we can conclude \eqref{eq:cor2} from an application of Young's inequality combined with Lemma \ref{lem:l1-bounds} and \eqref{eq:cor1}. 
\end{bew}

\subsection{Uniform bounds for the fluid component}
In preparation for obtaining uniform bounds on integrals involving $u_\epsi$, we will call for the following auxiliary result for ordinary differential equations as stated in \cite[Lemma 3.4]{SSW14}, where to we refer the reader for proof.
\begin{lemma}\label{lem:diffineq-lemma}
Let $T\in(1,\infty]$, $a>0$ and $b>0$. Suppose that $y:[0,T)\to[0,\infty)$ is absolutely continuous and such that
\begin{align*}
y'(t)+ay(t)\leq h(t)\quad\text{for a.e. }t\in(0,T),
\end{align*}
with some nonnegative $h\in\LSploc{1}{[0,T)}$ satisfying
\begin{align*}
\int_t^{t+1}h(s)\intd s\leq b\quad\text{for all }t\in[0,T-1).
\end{align*}
Then
\begin{align*}
y(t)\leq \max\left\{y(0)+b,\frac{b}{a}+2b\right\}\quad\text{for all }t\in(0,T).
\end{align*}
\end{lemma}

In quite standard manner (e.g. \cite[Lemmas 3.5 and 3.6]{win_globweak3d-AHPN16} or \cite[Lemma 4.3]{Wang17-globweak-ksns}) we can make use of the spatio-temporal estimate of $n_\epsi$ to obtain the following.

\begin{lemma}\label{lem:epsi-ind-est-u}
Let $m>\frac{4}{3}$ and assume that $n_0,c_0$ and $u_0$ comply with \eqref{IR}. Then there exists $C>0$ such that for all $\epsi\in(0,1)$ the solution $(n_\epsi,c_\epsi,u_\epsi)$ of \eqref{approxprob} satisfies
\begin{align*}
\intomega |u_\epsi|^2(\cdot,t)+\int_t^{t+1}\!\intomega|\nabla u_\epsi|^2\leq C
\end{align*}
for all $t>0$.
\end{lemma}

\begin{bew}
First, we note that in light of the Poincaré inequality and the embedding $W_{0,\sigma}^{1,2}(\Omega)\hookrightarrow\Lo[6]$ there exist $C_1>0$ and $C_2>0$ satisfying
\begin{align}\label{eq:epsi-ind-est-u-eq1}
\intomega |u_\epsi|^2(\cdot,t)\leq C_1\intomega |\nabla u_\epsi(\cdot,t)|^2\quad\text{for all }t>0\text{ and every }\epsi\in(0,1),
\end{align}
and
\begin{align}\label{eq:epsi-ind-est-u-eq2}
\|u_\epsi(\cdot,t)\|_{\Lo[6]}\leq C_2\intomega |\nabla u_\epsi|^2\quad\text{for all }t>0\text{ and all }\epsi\in(0,1).
\end{align}
Now, similar to the steps involving global existence (see \eqref{eq:globex-testing-eq4}), we test the third equation of \eqref{approxprob} by $u_\epsi$ and make use of integration by parts and Hölder's inequality to obtain 
\begin{align*}
\frac{1}{2}\frac{\intd}{\intd t}\intomega |u_\epsi|^2(\cdot,t)+\intomega|\nabla u_\epsi(\cdot,t)|^2\leq\|\nabla \phi\|_{\Lo[\infty]}\|u_\epsi(\cdot,t)\|_{\Lo[6]}\|n_\epsi(\cdot,t)\|_{\Lo[\nfrac{6}{5}]}
\end{align*}
for all $t>0$ and all $\epsi\in(0,1)$. Herein, we employ Young's inequality, \eqref{eq:epsi-ind-est-u-eq2} and \eqref{phidef} to find $C_3>0$ such that
\begin{align}\label{eq:epsi-ind-est-u-eq3}
\frac{1}{2}\frac{\intd}{\intd t}\intomega |u_\epsi|^2(\cdot,t)+\frac{1}{2}\intomega|\nabla u_\epsi(\cdot,t)|^2\leq C_3\|n_\epsi(\cdot,t)\|_{\Lo[\nfrac{6}{5}]}^2
\end{align}
holds for all $t>0$ and all $\epsi\in(0,1)$. Recalling that, by Lemma \ref{lem:st-bound-nepsi}, there exists $C_4>0$ satisfying $\int_t^{t+1}\|n_\epsi(\cdot,t)\|_{\Lo[\nfrac{6}{5}]}^2\leq\int_t^{t+1}\|n_\epsi(\cdot,t)+\epsi\|_{\Lo[\nfrac{6}{5}]}^2\leq C_4$ for all $t>0$ and all $\epsi\in(0,1)$, and estimating the gradient term on the left by means of \eqref{eq:epsi-ind-est-u-eq1}, an application of Lemma \ref{lem:diffineq-lemma} entails
\begin{align*}
\intomega |u_\epsi|^2(\cdot,t)\leq C_5\quad\text{for all }t>0\text{ and }\epsi\in(0,1),
\end{align*}
with some $C_5>0$. Returning to \eqref{eq:epsi-ind-est-u-eq3}, we integrate with respect to time to find that
\begin{align*}
\int_t^{t+1}\!\intomega|\nabla u_\epsi|^2\leq 2C_5+2C_3C_4\quad\text{for all }t>0,
\end{align*}
which concludes the proof.
\end{bew}

\subsection{Regularity estimates for the time derivatives}\label{sec42:timreg}
Obtaining information on the regularity of the time derivatives of our solution components is the next necessary step in preparing an Aubin--Lions type argument.
\begin{lemma}\label{lem:time-reg-n-and-c}
Let $m>\frac{4}{3}$ and suppose that $n_0,c_0$ and $u_0$ fulfill \eqref{IR}. For every $T>0$ there exists $C(T)>0$ such that for any $\epsi\in(0,1)$ the solution $(n_\epsi,c_\epsi,u_\epsi)$ of \eqref{approxprob} satisfies
\begin{align*}
\big\|\partial_t\big((n_{\epsi}+\epsi)^{m-1}\big)\big\|_{\LSp{1}{(0,T);(W_0^{3,2}(\Omega))^*}}\leq C(T),
\end{align*}
and
\begin{align*}
\|c_{\epsi t}\|_{\LSp{1}{(0,T);(W_0^{3,2}(\Omega))^*}}\leq C(T).
\end{align*}
\end{lemma}

\begin{bew}
Given $T>0$ we note that due to the continuous embedding of $\W[3,2]$ into  $\W[1,\infty]$ we can pick $C_1>0$ such that
\begin{align*}
\|\varphi\|_{\LSp{\infty}{(0,T);\W[1,\infty]}}\leq C_1\|\varphi\|_{\LSp{\infty}{(0,T);W_0^{3,2}(\Omega)}}\quad\text{for all }\varphi\in\LSp{\infty}{(0,T);W_0^{3,2}(\Omega)},
\end{align*}
Hence, for fixed $\varphi\in C_0^\infty(\Omega)$ with $\|\varphi\|_{W_0^{3,2}(\Omega)}\leq 1$ we can employ the Cauchy--Schwarz inequality to obtain
\begin{align}\label{eq:time-der-reg-nepsi-ineq}
\quad&\Big\vert\intomega \partial_t\big((n_{\epsi}+\epsi)^{m-1}\big)\varphi\Big\vert\nonumber\\
\leq\ &\frac{m|2-m|}{m-1}\intomega \big|\nabla(n_\epsi+\epsi)^{m-1}\big|^2|\varphi|+m\intomega (n_\epsi+\epsi)^{m-1}\big|\nabla (n_\epsi+\epsi)^{m-1}\cdot\!\nabla\varphi\big|\nonumber\\
 \ &+(m-1)|2-m|\!\intomega\frac{n_\epsi(n_\epsi+\epsi)^{m-3}}{(1+\epsi n_\epsi)^3}|\nabla n_\epsi\cdot\nabla c_\epsi|\,|\varphi|+(m-1)\!\intomega\frac{n_\epsi(n_\epsi+\epsi)^{m-2}}{(1+\epsi n_\epsi)^3}|\nabla c_\epsi\cdot\!\nabla\varphi|\nonumber\\ 
 \ &\quad+\intomega \big|u_\epsi\cdot\nabla (n_\epsi+\epsi)^{m-1}\big|\,|\varphi|\nonumber\\
\leq\ & C_1\frac{m|2-m|}{m-1}\intomega \big|\nabla(n_\epsi+\epsi)^{m-1}\big|^2+mC_1\Big(\intomega(n_\epsi+\epsi)^{2(m-1)}\Big)^\frac{1}{2}\Big(\intomega\big|\nabla (n_\epsi+\epsi)^{m-1}\big|^2\Big)^{\frac{1}{2}}\nonumber\\
\ &+C_1|2-m|\Big(\intomega\big|\nabla(n_\epsi+\epsi)^{m-1}\big|^2\Big)^\frac{1}{2}\Big(\intomega|\nabla c_\epsi|^2\Big)^\frac{1}{2}+C_1(m-1)\Big(\intomega (n_\epsi+\epsi)^{2(m-1)}\Big)^\frac{1}{2}\Big(\intomega|\nabla c_\epsi|^2\Big)^\frac{1}{2}\nonumber\\
&\quad+C_1\Big(\intomega |u_\epsi|^2\Big)^\frac{1}{2}\Big(\intomega\big|(\nabla n_\epsi+\epsi)^{m-1}\big|^2\Big)^\frac{1}{2}\quad\text{on }(0,T)\text{ for all }\epsi\in(0,1),
\end{align}
where we used the basic facts that $s\leq s+\epsi$ and $\frac{1}{(1+\epsi s)^3}\leq 1$ hold for all $s,\epsi\geq0$. Now, due to $2(m-1)<\frac{6m-4}{3}$, we infer from Young's inequality that
\begin{align*}
\intomega (n_\epsi+\epsi)^{2(m-1)}(\cdot,t)\leq \intomega (n_\epsi+\epsi)^{\frac{6m-4}{3}}(\cdot,t)+|\Omega|\quad\text{for all }t>0\text{ and all }\epsi\in(0,1),
\end{align*}
and hence, employing Young's inequality multiple times in \eqref{eq:time-der-reg-nepsi-ineq} and integrating with respect to time provides $C_2>0$ satisfying
\begin{align*}
\intoT\Big\vert\intomega \partial_t\big((n_{\epsi}+\epsi)^{m-1}\big)\varphi\Big\vert\leq C_2\intoTomega\big|\nabla(n_\epsi+\epsi)^{m-1}\big|^2&+C_2\intoTomega(n_\epsi+\epsi)^{\frac{6m-4}{3}}\\&+C_2\intoTomega|\nabla c_\epsi|^2+\intoTomega |u_\epsi|^2+C_2
\end{align*}
for all $\epsi\in(0,1)$ and all $\varphi\in\LSp{\infty}{(0,T);W_0^{3,2}(\Omega)}$ with $\|\varphi\|_{\LSp{\infty}{(0,T);W_0^{3,2}(\Omega)}}\leq1$. In consideration of  Lemma \ref{lem:bounds}, \eqref{eq:st-bound-specialcases} and Lemma \ref{lem:epsi-ind-est-u}, this entails the existence of $C_3(T)>0$ such that
\begin{align*}
\intoT\Big\vert\intomega \partial_t\big((n_\epsi+\epsi)^{m-1}\big)\varphi\Big\vert&\leq C_3(T)\quad\text{for all }\varphi\in\LSp{\infty}{(0,T);W_0^{3,2}(\Omega)}\text{ with }\|\varphi\|_{\LSp{\infty}{(0,T);W_0^{3,2}(\Omega)}}\leq1.
\end{align*}
In particular $\big\|\partial_t\big(n_\epsi+\epsi)^{m-1}\big)\big\|_{\LSp{1}{(0,T);(W_0^{3,2}(\Omega))^*}}\leq C_3(T)$, which we wanted to show. For the norm involving $c_{\epsi t}$ we work along similar lines, noticing that for fixed $\varphi$ as before we have
\begin{align}\label{eq:time-der-reg-cepsi}
\Big\vert\intomega c_{\epsi t}\varphi\Big\vert&\leq \intomega|\nabla c_\epsi\cdot\nabla \varphi|+\intomega c_\epsi\varphi+\intomega n_\epsi\varphi+\intomega c_\epsi|u_\epsi\cdot\nabla\varphi|\nonumber\\
&\leq C_1\intomega|\nabla c_\epsi|^2+C_1\intomega c_\epsi+C_1\intomega n_\epsi+\frac{C_1}{2}\intomega |u_\epsi|^2+\frac{C_1}{2}\intomega c_\epsi^2+C_4
\end{align}
holds with some $C_4>0$ on $(0,T)$ for all $\epsi\in(0,1)$. Thus, we conclude that the bounds contained in Lemma \ref{lem:l1-bounds}, Lemma \ref{lem:bounds} and Lemma \ref{lem:epsi-ind-est-u} immediately imply the the assertion upon integration of \eqref{eq:time-der-reg-cepsi} with respect to time.
\end{bew}

Relying on similar arguments one can also easily obtain a corresponding result for the fluid component.

\begin{lemma}\label{lem:time-reg-u}
Let $m>\frac{4}{3}$ and suppose that $n_0,c_0$ and $u_0$ fulfill \eqref{IR}. For every $T>0$ there exists $C(T)>0$ such that for any $\epsi\in(0,1)$ the solution $(n_\epsi,c_\epsi,u_\epsi)$ of \eqref{approxprob} satisfies
\begin{align}\label{eq:time-reg-u}
\int_0^T\|u_{\epsi t}\|^{\frac43}_{(W_{0,\sigma}^{1,2}(\Omega))^*}\leq C(T).
\end{align}
\end{lemma}

\begin{bew}
Following the reasoning of \cite[Lemma 5.5]{Wang17-globweak-ksns}, we fix an arbitrary $\psi\in C_0^\infty(\Omega)$ with $\nabla \cdot\psi\equiv0$ in $\Omega$ and make use of the third equation in \eqref{approxprob} and Hölder's inequality, to find that
\begin{align}\label{eq:time-reg-u-1}
\Big|\intomega u_{\epsi t}\cdot\psi\Big|\leq\|\nabla u_\epsi\|_{\Lo[2]}\|\nabla\psi\|_{\Lo[2]}\nonumber+\|Y_\epsi u_\epsi\|_{\Lo[6]}&\|u_\epsi\|_{\Lo[3]}\|\nabla\psi\|_{\Lo[2]}\\&+\|\nabla\phi\|_{\Lo[\infty]}\|n_\epsi\|_{\Lo[\frac{6}{5}]}\|\psi\|_{\Lo[6]}
\end{align}
holds  on $(0,\infty)$ for all $\epsi\in(0,1)$. To further estimate the norm of the $Y_\epsi u_\epsi$, we make use of the embedding $W_{0,\sigma}^{1,2}(\Omega)\hookrightarrow\Lo[6]$, as well as the facts that $Y_\epsi$ commutes with $A^{\frac{1}{2}}$ and is nonexpansive on $L_\sigma^2(\Omega)$ to obtain
\begin{align*}
\|Y_\epsi u_\epsi(\cdot,t)\|_{\Lo[6]}&\leq\|\nabla Y_\epsi u_\epsi(\cdot,t)\|_{\Lo[2]}=\|A^\frac{1}{2}Y_\epsi u_\epsi(\cdot,t)\|_{\Lo[2]}\\&=\|Y_\epsi A^{\frac{1}{2}}u_\epsi(\cdot,t)\|_{\Lo[2]}\leq\|A^\frac{1}{2}u_\epsi(\cdot,t)\|_{\Lo[2]}=\|\nabla u_\epsi(\cdot,t)\|_{\Lo[2]}
\end{align*}
for all $t>0$ and any $\epsi\in(0,1)$. Combination of this with \eqref{eq:time-reg-u-1}, \eqref{phidef} and the boundedness of $\psi$ and its derivative, entails the existence of $C_1>0$ such that
\begin{align}\label{eq:time-reg-u-2}
\|u_{\epsi t}(\cdot,t)\|_{(W_{0,\sigma}^{1,2}(\Omega))^*}^\frac{4}{3}\leq C_1\Big(\|\nabla u_\epsi(\cdot,t)\|_{\Lo[2]}^\frac{4}{3}+\|\nabla u_\epsi(\cdot,t)\|_{\Lo[2]}^\frac{4}{3}\|u_\epsi(\cdot,t)\|_{\Lo[3]}^\frac{4}{3}+\|n_\epsi(\cdot,t)\|_{\Lo[\frac{6}{5}]}^\frac{4}{3}\Big)
\end{align}
for all $t>0$ and all $\epsi\in(0,1)$. Now, in light of the Young and Gagliardo--Nirenberg inequalities we have
\begin{align*}
\|n_\epsi(\cdot,t)\|_{\Lo[\frac{6}{5}]}^\frac{4}{3}\leq \|n_\epsi(\cdot,t)\|_{\Lo[\frac{6}{5}]}^2+C_2\quad\text{ and }\quad
\|u_\epsi(\cdot,t)\|_{\Lo[3]}^\frac{4}{3}\leq C_3\|\nabla u_\epsi(\cdot,t)\|_{\Lo[2]}^\frac{2}{3}\|u_\epsi(\cdot,t)\|_{\Lo[2]}^\frac{2}{3}
\end{align*}
for all $\epsi\in(0,1)$ on $(0,\infty)$, with some $C_2>0$ and $C_3>0$. Hence, plugging these two estimates into \eqref{eq:time-reg-u-2} and integrating with respect to time we obtain
\begin{align*}
\intoT\|u_{\epsi t}\|_{(W_{0,\sigma}^{1,2}(\Omega))^*}^\frac{4}{3}\leq C_1\intoT\|\nabla u_\epsi\|_{\Lo[2]}^\frac{4}{3}+C_1C_3\intoT\|\nabla u_\epsi\|_{\Lo[2]}^2\|u_\epsi\|_{\Lo[2]}^\frac{2}{3}+C_1\intoT\|n_\epsi\|_{\Lo[\frac{6}{5}]}^2+C_1C_2T
\end{align*}
for all $T>0$ and all $\epsi\in(0,1)$, and thus \eqref{eq:time-reg-u} is an evident consequence of the bounds featured in Lemma \ref{lem:st-bound-nepsi} and Lemma \ref{lem:epsi-ind-est-u}.
\end{bew}

\setcounter{equation}{0} 
\section{Existence of limit functions}\label{sec5:convergence}
With the uniform bounds from Lemma \ref{lem:bounds}, Lemma \ref{lem:epsi-ind-est-u}, Lemma \ref{lem:time-reg-n-and-c} and Lemma \ref{lem:time-reg-u} we are now in the position to obtain limit functions $n,c$ and $u$, which at least fulfill the regularity assumptions required in Definition \ref{def:very_weak_sol}.

\begin{lemma}\label{lem:convergences}
Let $m>\frac{4}{3}$ and suppose that $n_0,c_0,u_0$ comply with \eqref{IR}.  Then there exist a sequence $(\epsi_j)_{j\in\N}\subset(0,1)$ with $\epsi_j\searrow0$ as $j\to\infty$ and functions
\begin{align*}
n&\in\LSploc{2m-\frac{4}{3}}{\bomega\times[0,\infty)}\quad\text{with}\quad\nabla n^{m-1}\in\LSploc{2}{\bomega\times[0,\infty)},\\
c&\in\LSploc{2}{[0,\infty);\W[1,2]},\\
u&\in\LSploc{2}{[0,\infty);W^{1,2}_{0,\sigma}(\Omega)},
\end{align*}
such that the solutions $(n_\epsi,c_\epsi,u_\epsi)$ of \eqref{approxprob} satisfy
\begin{alignat}{2}
(n_\epsi+\epsi)^{m-1}&\to n^{m-1}&&\text{in }\LSploc{2}{\bomega\times[0,\infty)}\text{ and a.e. in }\Omega\times(0,\infty),\label{eq:conv-n-m-1-ae}\\
\nabla (n_\epsi+\epsi)^{m-1}&\wto \nabla n^{m-1}&&\text{in }\LSploc{2}{\bomega\times[0,\infty)},\label{eq:conv-nab-n-m-1}\\
n_\epsi+\epsi&\wto n\quad\qquad&&\text{in }\LSploc{2m-\frac{4}{3}}{\bomega\times[0,\infty)},\label{eq:conv-nw}\\
n_\epsi+\epsi\to n\quad\text{and}\quad n_\epsi&\to n&&\text{in }\LSploc{p}{\bomega\times[0,\infty)}\text{ for any }p\in[1,2m-\tfrac43),\label{eq:conv-n-strong}\\
c_\epsi&\to c&&\text{in }\LSploc{2}{\bomega\times[0,\infty)}\text{ and a.e in }\Omega\times(0,\infty),\label{eq:conv-c}\\
\nabla c_\epsi&\wto \nabla c&&\text{in }\LSploc{2}{\bomega\times[0,\infty)},\label{eq:conv-nab-c}
\intertext{as well as}
u_\epsi&\to u&&\text{in }\LSploc{2}{\bomega\times[0,\infty)},\label{eq:conv-u-l2}\\
\nabla u_\epsi&\wto \nabla u&&\text{in }\LSploc{2}{\bomega\times[0,\infty)},\label{eq:conv-nab-u}\\
Y_\epsi u_\epsi&\to u&&\text{in }\LSploc{2}{\bomega\times[0,\infty)}\label{eq:conv-Yu}
\end{alignat}
as $\epsi=\epsi_j\searrow0$, and such that $n\geq0$, $c\geq0$ a.e. in $\Omega\times(0,\infty)$. If, moreover, $m\in(\tfrac43,2)$, then there exists a further subsequence $(\epsi_{j_k})_{k\in\N}\subset(0,1)$ such that $(n_\epsi,c_\epsi,u_\epsi)$ also satisfy
\begin{alignat}{2}
(n_\epsi+1)^{m-1}&\to (n+1)^{m-1}&&\text{in }\LSploc{2}{\bomega\times[0,\infty)},\label{eq:conv-n+1-m-1-ae}\\
\nabla (n_\epsi+1)^{m-1}&\wto \nabla (n+1)^{m-1}&&\text{in }\LSploc{2}{\bomega\times[0,\infty)},\label{eq:conv-nab-n+1-m-1}\\
(n_\epsi+1)^{\frac{m-3}{2}}(n_\epsi+\epsi)^{\frac{m-1}{2}}\nabla n_\epsi&\wto (n+1)^{\frac{m-3}{2}}n^{\frac{m-1}{2}} \nabla n\qquad\quad&&\text{in }\LSploc{2}{\bomega\times[0,\infty)},\qquad\label{eq:conv-nab-n+1+eps-m-1}
\end{alignat}
as $\epsi=\epsi_{j_k}\searrow0$.
\end{lemma}

\begin{bew}
Since $2(m-1)<2m-\frac{4}{3}$, the bounds featured in Lemma \ref{lem:st-bound-nepsi}, Lemma \ref{lem:bounds} and Lemma \ref{lem:time-reg-n-and-c} imply that
\begin{align*}
\big\{(n_\epsi+\epsi)^{m-1}\big\}_{\epsi\in(0,1)}\quad\text{is bounded in}\ \LSploc{2}{[0,\infty);\W[1,2]}
\end{align*}
and that
\begin{align*}
\big\{\partial_t (n_\epsi+\epsi)^{m-1}\big\}_{\epsi\in(0,1)}\quad\text{is bounded in }\ \LSploc{1}{[0,\infty);(W_0^{3,2}(\Omega))^*}.
\end{align*}
Hence, an Aubin--Lions type lemma (e.g. \cite[Corollary 8.4]{Sim87}) provides the existence of $(\epsi_j)_{j\in\N}\subset(0,1)$ satisfying $\epsi_j\searrow 0$ as $j\to\infty$ such that \eqref{eq:conv-n-m-1-ae} holds. The weak convergences stated in \eqref{eq:conv-nab-n-m-1} and \eqref{eq:conv-nw} are immediate consequences of the spatio-temporal bounds contained in Lemma \ref{lem:bounds} and Lemma \ref{lem:st-bound-nepsi}, respectively, upon extraction of a further (non-relabeled) subsequence, whereas the improvement to strong convergence obtained in the first part of \eqref{eq:conv-n-strong} follows from an application of the Vitali convergence theorem while relying on the a.e. convergence of $n_\epsi+\epsi$ entailed by \eqref{eq:conv-n-m-1-ae} and the equi-integrability property of $\{(n_{\epsi_j}+\epsi_j)^p\}_{j\in\N}$ for $p<2m-\frac{4}{3}$ contained in \eqref{eq:st-bound-specialcases}. The second part of \eqref{eq:conv-n-strong} then is an obvious consequence of the uniform convergence of $\epsi_j$ to zero. In a similar fashion, we make use of Lemmas \ref{lem:bounds} and \ref{lem:time-reg-n-and-c} in combination with the Aubin--Lions lemma to find that and \eqref{eq:conv-c} and \eqref{eq:conv-nab-c} hold. Applying the same arguments to bounds on $u_\epsi$ in $\LSploc{2}{[0,\infty);W_\sigma^{1,2}(\Omega)}$ and of $u_{\epsi t}$ in $\LSploc{\frac43}{[0,\infty),(W_{0,\sigma}^{1,2}(\Omega))^*}$, as implied by Lemma \ref{lem:epsi-ind-est-u} and Lemma \ref{lem:time-reg-u}, also proves \eqref{eq:conv-u-l2} and \eqref{eq:conv-nab-u}. 
Finally, relying on arguments as in e.g. \cite[Lemma 4.1]{win_globweak3d-AHPN16}, we make use of the properties of the Yosida approximation (see \cite[II.3.4]{sohr}) to find that for all $\varphi\in L_{\sigma}^2(\Omega)$ we have $\|Y_{\epsi}\varphi\|_{\Lo[2]}\leq\|\varphi\|_{\Lo[2]}$ and $Y_\epsi\varphi\to\varphi$ in $\Lo[2]$ as $\epsi\searrow0$ to conclude from \eqref{eq:conv-u-l2} that 
\begin{align*}
\|Y_{\epsi_j} u_{\epsi_j}(\cdot,t)-u(\cdot,t)\|_{\Lo[2]}\leq \|u_{\epsi_j}(\cdot,t)-u(\cdot,t)\|_{\Lo[2]}+\|Y_{\epsi_j}u(\cdot,t)-u(\cdot,t)\|_{\Lo[2]}\to 0
\end{align*}
for a.e. $t>0$ as $\epsi_j\searrow0$. Since moreover, $\|Y_{\epsi_j} u_{\epsi_j}(\cdot,t)-u(\cdot,t)\|_{\Lo[2]}^2\leq 4\sup_{\epsi\in(0,1)}\|u_\epsi\|_{\LSp{2}{\Omega\times(0,\infty)}}^2$ for a.e $t>0$, an application of the dominated convergence theorem implies \eqref{eq:conv-Yu} in light of Lemma \ref{lem:epsi-ind-est-u}. The remaining convergence properties for $m\in(\tfrac43,2)$ follow in a similar fashion from the bounds contained in Corollary \ref{cor:1}, where for the convergence statement in \eqref{eq:conv-n+1-m-1-ae} we once more rely on Vitali's theorem.
\end{bew}

\setcounter{equation}{0} 
\section{Solution properties of the limit functions}\label{sec6:sol-prop}

\subsection{Weak solution properties of \texorpdfstring{$c$ and $u$}{c and u}}\label{sec61:weak-sol-cu}
As an immediate consequence of the convergences presented in Lemma \ref{lem:convergences} we also obtain the following.
\begin{lemma}\label{lem:sol-prop-c-u}
Let $m>\frac{4}{3}$ and assume that $n_0,c_0$ and $u_0$ comply with \eqref{IR}. Furthermore, let $n,c,u$ denote the limit functions provided by Lemma \ref{lem:convergences}. Then 
\begin{align}\label{eq:mass-ineq}
\intomega n(\cdot,t)= \intomega n_0\quad\text{for a.e. }t>0,
\end{align}
and $c$ and $u$ satisfy the weak solution properties \eqref{eq:very_weak_sol1} and \eqref{eq:very_weak_sol2}, respectively, of Definition \ref{def:very_weak_sol}.
\end{lemma}

\begin{bew}
The mass equality \eqref{eq:mass-ineq} is an immediate consequence of the strong convergence statement \eqref{eq:conv-n-strong} and Lemma \ref{lem:l1-bounds}. Testing the second equation of \eqref{approxprob} with an arbitrary compactly in $\bomega\times[0,\infty)$ supported test function $\varphi\in\LSp{\infty}{\Omega\times(0,\infty)}\cap\LSp{2}{(0,\infty);\W[1,2]}$ with $\varphi_t\in\LSp{2}{\Omega\times(0,\infty)}$ we obtain
\begin{align*}
-\intinfomega c_\epsi\varphi_t-\intomega c_0\varphi(\cdot,0)=-\intinfomega\nabla c_\epsi\cdot\nabla\varphi-\intinfomega c_\epsi\varphi+\intinfomega n_\epsi\varphi+\intinfomega c_\epsi (u_\epsi\cdot\nabla\varphi)
\end{align*}
for all $\epsi\in(0,1)$. Recalling that the convergence properties contained in Lemma \ref{lem:convergences}, in particular \eqref{eq:conv-c}, \eqref{eq:conv-nab-c}, \eqref{eq:conv-n-strong} and \eqref{eq:conv-u-l2}, are clearly sufficient to pass to the limit in all of the integrals, we conclude \eqref{eq:very_weak_sol1}. Similarly, testing the third equation of \eqref{approxprob} by an arbitrary $\psi\in C_0^\infty(\Omega\times[0,\infty);\R^3)$ with $\nabla\cdot\psi\equiv0$ in $\Omega\times(0,\infty)$ we find that
\begin{align*}
-\intinfomega u_\epsi\psi_t-\intomega u_0\psi(\cdot,0)=-\intinfomega\nabla u_\epsi\cdot\nabla\psi+\intinfomega (Y_\epsi u_\epsi\otimes u_\epsi)\cdot\nabla\psi+\intinfomega n_\epsi(\nabla\phi\cdot\psi)
\end{align*}
holds for all $\epsi\in(0,1)$. Relying on the convergence properties \eqref{eq:conv-u-l2}, \eqref{eq:conv-nab-u}, \eqref{eq:conv-Yu} and \eqref{eq:conv-n-strong} obtained in Lemma \ref{lem:convergences}, where specifically \eqref{eq:conv-u-l2} and \eqref{eq:conv-Yu} also entail that $Y_\epsi u_\epsi\otimes u_\epsi\to u\otimes u$ in $\LSploc{2}{\bomega\times[0,\infty)}$, we can pass to the limit in all integrals and infer that \eqref{eq:very_weak_sol2} is valid.
\end{bew}

\subsection{Weak solution property of \texorpdfstring{$n$ for $m>\frac53$}{n for m>5/3}}\label{sec62:weak-sol-n}
Recalling that by \eqref{eq:conv-n-strong} we have $n_\epsi\to n$ in $\LSploc{2}{\bomega\times[0,\infty)}$, whenever $m$ is sufficiently large, we observe that hence the weak convergence results obtained for $\nabla(n_\epsi+\epsi)^{m-1}$ and $\nabla c_\epsi$ are already sufficient to show that \eqref{eq:very_weak_supersol} holds for $\Phi(s)\equiv s$ with equality and that hence the solution is in fact a global weak solution in the standard sense.
\begin{lemma}\label{lem:weak-sol-5over3}
Let $m>\frac{5}{3}$, suppose that $n_0,c_0$ and $u_0$ comply with \eqref{IR}, and let $n,c,u$ denote the limit functions obtained in Lemma \ref{lem:convergences}. Then $n\in\LSploc{2}{\bomega\times[0,\infty)}$ and for any $\varphi\in C_0^\infty\big(\bomega\times[0,\infty)\big)$ the weak solution property \eqref{eq:weak-sol-n} is satisfied.
\end{lemma}

\begin{bew}
Testing the first equation of \eqref{approxprob} by $\varphi\in C_0^\infty\big(\bomega\times[0,\infty)\big)$ we find that $n_\epsi$ satisfies
\begin{align}\label{eq:proof-weak-sol-n-eq1}
-\intinfomega n_\epsi\varphi_t-\intomega n_0\varphi(\cdot,0)=-\frac{m}{m-1}&\intinfomega (n_\epsi+\epsi)\big(\nabla(n_\epsi+\epsi)^{m-1}\cdot\nabla\varphi\big)\\+&\intinfomega\frac{n_\epsi}{(1+\epsi n_\epsi)^3}(\nabla c_\epsi\cdot\nabla\varphi)+\intinfomega n_\epsi(u_\epsi\cdot\nabla\varphi)\nonumber
\end{align}
for all $\epsi\in(0,1)$. Since $m>\frac{5}{3}$ implies $2m-\frac{4}{3}>2$, we obtain from \eqref{eq:conv-n-strong} that 
\begin{align*}
(n_\epsi+\epsi)\to n\quad\text{and}\quad n_\epsi\to n\quad\text{in}\quad \LSploc{2}{\bomega\times[0,\infty)}\quad\text{as }\epsi=\epsi_j\searrow0.
\end{align*}
Making additional use of \eqref{eq:conv-nab-n-m-1} and \eqref{eq:conv-u-l2}, we thus have
\begin{align*}
-\frac{m}{m-1}\intinfomega(n_\epsi+\epsi)\big(\nabla(&n_\epsi+\epsi)^{m-1}\cdot\nabla\varphi\big)\to -\frac{m}{m-1}\intinfomega n\big(\nabla n^{m-1}\cdot\nabla\varphi\big),\\
-\intinfomega n_\epsi\varphi_t\to-\intinfomega n\varphi_t,\quad&\text{and}\quad\intinfomega n_\epsi(u_\epsi\cdot\nabla\varphi)\to\intinfomega n(u\cdot\nabla\varphi),\quad\text{as }\epsi=\epsi_j\searrow0.
\end{align*}
To treat the remaining integral we note that since $\frac{1}{(1+\epsi n_\epsi)^3}\leq 1$ for all $\epsi\in(0,1)$ and $\frac{1}{(1+\epsi n_\epsi)^3}\to 1$ a.e. in $\Omega\times(0,\infty)$, we can employ a useful effect of the dominated convergence theorem (see \cite[Lemma 10.4]{win15_chemorot}) to find that
\begin{align*}
\frac{n_\epsi}{(1+\epsi n_\epsi)^3}\to n\quad\text{in}\quad\LSploc{2}{\bomega\times[0,\infty)}\quad\text{as }\epsi=\epsi_j\searrow0,
\end{align*}
which combined with \eqref{eq:conv-nab-c} shows
\begin{align*}
\intinfomega \frac{n_\epsi}{(1+\epsi n_\epsi)^3}(\nabla c_\epsi\cdot\nabla \varphi)\to\intinfomega n(\nabla c\cdot\nabla\varphi)\quad\text{as }\epsi=\epsi_{j}\searrow0.
\end{align*}
In conclusion, we may take $\epsi=\epsi_{j}\searrow0$ in \eqref{eq:proof-weak-sol-n-eq1} to find that \eqref{eq:weak-sol-n} is valid.
\end{bew}

The combination of three of our previous results now immediately establishes Theorem \ref{theo:1}.

\begin{proof}[\textbf{Proof of Theorem \ref{theo:1}}:]
We can merge the results of Lemma \ref{lem:sol-prop-c-u} and Lemma \ref{lem:weak-sol-5over3} with the regularity properties contained in Lemma \ref{lem:convergences} to immediately arrive at the conclusion.
\end{proof}

\subsection{Very weak solution property of \texorpdfstring{$n$}{n} in the case of \texorpdfstring{$m>\frac{4}{3}$}{m>4/3}}\label{sec63:very-weak-sol}
Since the strong convergence of $n_\epsi$ in $\LSploc{2}{\bomega\times[0,\infty)}$ is heavily reliant on the fact that $m>\frac{5}{3}$, we cannot expect global weak solutions for $m\leq\frac{5}{3}$. Having in mind global very weak solutions as defined in Definition \ref{def:very_weak_sol} instead, we find that our current precompactness properties are insufficient to treat some of the terms arising in \eqref{eq:very_weak_supersol}. In particular, since we only have a weak convergence for $\nabla n_\epsi^{m-1}$ in $\LSploc{2}{\bomega\times[0,\infty)}$ at hand, we have to improve our convergence result for $\nabla c_\epsi$ in order to treat the mixed derivative term. As a preparatory result, we state the following Lemma, which has been proven in \cite[Lemma 7.1]{Wang17-globweak-ksns} in a closely related setting.
\begin{lemma}\label{lem:c2-ineq}
Let $m>\frac{4}{3}$ and assume that $n_0,c_0$ and $u_0$ comply with \eqref{IR}. Then there exists a null set $N\subset(0,\infty)$ such that the functions $n,c$ and $u$ obtained in Lemma \ref{lem:convergences} satisfy
\begin{align}\label{eq:c2-ineq}
\frac{1}{2}\intomega c^2(\cdot,T)-\frac{1}{2}\intomega c_0^2+\intoTomega|\nabla c|^2\geq -\intoTomega c^2+\intoTomega nc\quad\text{for all }T\in(0,\infty)\setminus N.
\end{align}
\end{lemma}

\begin{bew}
The arguments employed to prove the asserted inequality can in detail be found in \cite[Lemma 7.1]{Wang17-globweak-ksns}, which  adapts the reasoning found in \cite[Lemma 8.1]{win15_chemorot} to the signal production setting while also including fluid terms. As the proof is quite technical and in essence unchanged, we will refrain from a detailed reconstruction of the proof and only sketch the main steps. For more details the reader is referred to \cite{Wang17-globweak-ksns,win15_chemorot}.

Since \eqref{eq:conv-c} shows that $z(t):=\intomega c^2(\cdot,t)$, $t>0$, satisfies $z\in\LSploc{1}{[0,\infty)}$, we can find a null set $N\subset(0,\infty)$ such that each $T\in(0,\infty)\setminus N$ is a Lebesgue point of $z$, which in turn shows that
\begin{align*}
\frac{1}{\delta}\int_{T}^{T+\delta}\!\intomega c^2(\cdot,t)\to\intomega c^2(\cdot,T)\quad\text{for all }T\in(0,\infty)\setminus N\ \text{ as }\ \delta\searrow 0.
\end{align*} 
To prepare a special test function to use in the second equation of \eqref{CTnonlindiff}, for given $T\in(0,\infty)\setminus N$ and $\delta\in(0,1)$ and $r\in(0,1)$ we let
\begin{align*}
\zeta_\delta(t):=\begin{cases}
1,\quad& t\in[0,T],\\
\frac{T+\delta-t}{\delta},&t\in(T,T+\delta),\\
0,& t\geq T+\delta,
\end{cases}\quad\text{and}\quad \psi_r(s):=\frac{s}{1+rs},
\end{align*}
as well as
\begin{align*}
\tilde{c}_k(x,t):=\begin{cases}
c(x,t),\quad&(x,t)\in\Omega\times(0,\infty),\\
c_{0k}(x),&(x,t)\in\Omega\times(-1,0],
\end{cases}
\end{align*}
for $k\in\N$, where the nonnegative sequence $(c_{0k})_{k\in\N}\subset \CSp{1}{\bomega}$ is chosen such that $c_{0k}\to c_0$ in $\Lo[2]$ as $k\to\infty$. Denoting by
\begin{align*}
\big(A_h\psi_r(\tilde{c}_k)\big)(x,t):=\frac{1}{h}\int_{t-h}^t\psi_r(\tilde{c}_k)(x,s)\intd s,\quad(x,t)\in\Omega\times[0,\infty),
\end{align*}
the temporal average and with $h\in(0,1)$ letting
\begin{align*}
\varphi(x,t):=\varphi_{\delta,k,h,r}(x,t):=\zeta_\delta(t)\cdot\big(A_h\psi_r(\tilde{c}_k)\big)(x,t),\quad(x,t)\in\Omega\times[0,\infty),
\end{align*}
we can check that $\varphi$ is of class $\LSploc{\infty}{\bomega\times[0,\infty)}\cap\LSp{2}{(0,\infty);\W[1,2]}$, that $\varphi$ has compact support in $\bomega\times[0,T+1]$ and that $\varphi\in\LSp{2}{\Omega\times(0,\infty)}$ and hence, $\varphi$ is an admissible test function for \eqref{eq:very_weak_sol1}. Inserting $\varphi$ into \eqref{eq:very_weak_sol1} we obtain
\begin{align}\label{eq:strong-conv-ineq-proof1}
-\intinfomega \zeta_\delta'(t)\cdot&\big(A_h\psi_r(\tilde{c}_k)\big)(x,t)\cdot c(x,t)\intd x\intd t-\!\intinfomega\frac{\zeta_\delta(t)}{h}\!\cdot\!\big[\psi_r(\tilde{c}_k)(x,t)-\psi_r(\tilde{c}_k)(x,t-h)\big]\!\cdot\! c(x,t)\intd x\intd t\nonumber \\-\intomega c_0(x)\cdot&\big(A_h\psi_r(\tilde{c}_k)\big)(x,0)\intd x=-\intinfomega\nabla c(x,t)\cdot\zeta_\delta(t)\cdot\nabla\big(A_h\psi_r(\tilde{c}_k)\big)(x,t)\intd x\intd t\nonumber\\-\intinfomega c(x,t)&\cdot\zeta_\delta(t)\cdot\big(A_h\psi_r(\tilde{c}_k)\big)(x,t)\intd x\intd t+\intinfomega n(x,t)\cdot\zeta_\delta(t)\cdot\big(A_h\psi_r(\tilde{c}_k)\big)(x,t)\intd x\intd t\nonumber\\
+\intinfomega c(x,t)&\cdot\zeta_\delta(t) \cdot u(x,t)\cdot\nabla \big(A_h\psi_r(\tilde{c}_k)\big)(x,t)\intd x\intd t.
\end{align}
Here, we note that $\psi_r(\tilde{c}_k)\in\LSp{\infty}{\Omega\times(0,T+1)}$, that the fact that $c_{0k}\in\CSp{1}{\bomega}$ implies $\nabla\psi_r(\tilde{c}_k)\in\LSp{2}{\Omega\times(0,T+1)}$, and that $\Psi_r(s):=\frac{rs-\ln(1+rs)}{r^2}$ is the primitive of $\psi_r(s)$ for any $s\geq0$. Hence, we can make use of known results for Steklov averages (see e.g. \cite[Lemma 10.2]{win15_chemorot}) to let $h\searrow0$ in \eqref{eq:strong-conv-ineq-proof1} and obtain
\begin{align}\label{eq:strong-conv-ineq-proof2}
\ &-\liminf_{h\to0}\intinfomega\frac{\zeta_\delta(t)}{h}\cdot\big[\psi_r(\tilde{c}_k(x,t)-\psi_r(\tilde{c}_k(x,t-h)\big]\cdot c(x,t)\intd x\intd t\nonumber\\=\ &-\intinfomega\zeta_\delta(t)\frac{|\nabla c(x,t)|^2}{(1+rc(x,t))^2}\intd x\intd t-\intinfomega\zeta_\delta(t)\frac{c^2(x,t)}{1+rc(x,t)}\intd x\intd t\\&+\intinfomega\zeta_\delta(t)\frac{n(x,t)c(x,t)}{1+rc(x,t)}\intd x\intd t+\intinfomega\zeta_\delta'(t)\frac{c^2(x,t)}{1+rc(x,t)}\intd x\intd t+\intomega\frac{c_0(x)c_{0k}(x)}{1+r c_{0k}(x)}\intd x.\nonumber
\end{align}
To estimate the remaining limit (compare (7.11)--(7.14) in \cite[Lemma 7.1]{Wang17-globweak-ksns}), we make use of the convexity of $\Psi_r$ implying
\begin{align*}
\Psi_r(\tilde{c}_k(x,t+h))-\Psi_r(\tilde{c}_k(x,t))\geq\psi_r(\tilde{c}_k)(x,t)\big(c(x,t+h)-c(x,t)\big)\text{ for a.e. }x\in\Omega\text{ and }t\in(0,T+1), 
\end{align*}
as well as the substitution $s=t+h$, Young's inequality and the definition of $\zeta_\delta$ to find that
\begin{align*}
\ &-\intinfomega\frac{\zeta_\delta(t)}{h}\cdot\big[\psi_r(\tilde{c}_k(x,t)-\psi_r(\tilde{c}_k(x,t-h)\big]\cdot c(x,t)\intd x\intd t\\ \leq\ &
\intinfomega\frac{\zeta_\delta(t+h)}{h}\big[\Psi_r(\tilde{c}_k(x,t+h))-\Psi_r(\tilde{c}_k(x,t))\big]\intd x\intd t +\frac{1}{2}\intomega\frac{c_{0k}^2(x)}{(1+r c_{0k}(x))^2}\intd x\\\ &\hspace*{3.6cm}+\frac{1}{2h}\int_0^h\!\intomega c^2(x,t)\intd x\intd t+\intinfomega\frac{\zeta_\delta(t+h)-\zeta_\delta(h)}{h}\psi_r(\tilde{c}_k(x,t))c(x,t)\intd x\intd t,
\end{align*}
which upon combination with \eqref{eq:strong-conv-ineq-proof2} shows that
\begin{align*}
\ &\intinfomega\zeta_\delta(t)|\nabla c(x,t)|^2\intd x\intd t+\frac{1}{2}\intomega c^2_{0k}(x)\intd x+\frac{1}{2}\intomega c_0^2(x)\intd x\\\geq\ &-\intinfomega\zeta_\delta(t)\frac{c^2(x,t)}{1+rc(x,t)}\intd x\intd t+\intinfomega\zeta_\delta(t)\frac{n(x,t)c(x,t)}{1+rc(x,t)}\intd x\intd t+\intomega\frac{c_0(x)c_{0k}(x)}{1+rc_{0k}(x)}\intd x\\&+\intinfomega\zeta_\delta'(t)\Psi_r(\tilde{c}_k(x,t))\intd x\intd t+\!\intomega\Psi_r(\tilde{c}_k(x,0))\intd x
\end{align*}
for all $k\in\N$ and $r\in(0,1)$. By means of the dominated convergence theorem, we may next let $r\searrow0$ and then $k\to\infty$ to arrive at
\begin{align*}
\ &\intinfomega\zeta_\delta(t)|\nabla c(x,t)|^2\intd x\intd t\intd x+\intinfomega\zeta_\delta(t)c^2(x,t)\intd x\intd t-\intinfomega\zeta_\delta(t)n(x,t)c(x,t)\intd x\intd t\\\geq\ \ &\frac{1}{2}\intomega c_0^2(x)\intd x-\frac{1}{2\delta}\int_T^{T+\delta}\!\intomega c^2(x,t)\intd x\intd t.
\end{align*}
Finally, recalling the Lebesgue point property of $T$ we make use of the dominated convergence theorem once more to take $\delta\searrow0$ and obtain \eqref{eq:c2-ineq}.
\end{bew}

The inequality from the previous lemma at hand, we can now make use of arguments previously employed in \cite[Lemma 4.4]{win15_stokesrot} and \cite[Lemma 7.2]{Wang17-globweak-ksns} to obtain the last missing convergence property we require in order to pass to the limit in the integrals appearing in the very weak $\Phi$--supersolution concept.

\begin{lemma}\label{lem:strong-conv-nabc-l2}
Let $m>\frac{4}{3}$ and assume that $n_0,c_0$ and $u_0$ comply with \eqref{IR}. Furthermore, denote by $(\epsi_j)_{j\in\N}$ and $n,c,u$ the sequence and limit functions provided by Lemma \ref{lem:convergences}. Then there exist a subsequence $(\epsi_{j_k})_{k\in\N}$ and a null set $N\subset(0,\infty)$ such that for each $T\in(0,\infty)\setminus N$ the classical solution $(n_\epsi,c_\epsi,u_\epsi)$ of \eqref{approxprob} satisfies
\begin{align*}
\nabla c_\epsi\to\nabla c\quad\text{in }\LSp{2}{\Omega\times(0,T)}\ \text{as}\ \epsi=\epsi_{j_k}\searrow0.
\end{align*}
\end{lemma}

\begin{bew}
To start, let us set $l:=2-\frac{4}{3m}$, which by the assumption $m>\frac{4}{3}$ satisfies $2>l>1$, as well as
\begin{align*}
\frac{6l}{6-l}=\frac{9m-6}{3m+1}\in(1,6(m-1)),\quad\text{and}\quad\frac{2l}{2-l}=3m-2=\frac{2\cdot\frac{6l}{6-l}\cdot(m-\frac{7}{6})}{\frac{6l}{6-l}-1},
\end{align*}
which in consequence of Lemma \ref{lem:st-bound-nepsi} implies that
\begin{align}\label{eq:strong-conv-nabc-l2-proof-eq1}
\int_t^{t+1}\|n_\epsi(\cdot,s)\|_{\Lo[\frac{6l}{6-l}]}^{\frac{2l}{2-l}}\intd s\leq C_1\quad\text{for all }t>0\text{ and all }\epsi\in(0,1),
\end{align}
with some $C_1>0$. Now, with $N_1\subset(0,\infty)$ denoting the null set obtained in Lemma \ref{lem:c2-ineq} we note that by Lemma \ref{lem:convergences} we can find another null set $N_2\supset N_1$ and a subsequence $(\epsi_{j_k})_{k\in\N}$ such that
\begin{align*}
\intomega c_\epsi^2(\cdot,T)\to\intomega c^2(\cdot,T)\quad\text{for all }T\in(0,\infty)\setminus N_2\quad\text{as}\quad \epsi=\epsi_{j_k}\searrow0.
\end{align*}
For any such $T\in(0,\infty)\setminus N_2$ we find by the Hölder and Young inequalities that
\begin{align*}
\intoTomega |n_\epsi c_\epsi|^l&\leq\intoT\|n_\epsi(\cdot,s)\|_{\Lo[\frac{6l}{6-l}]}^l\|c_\epsi(\cdot,s)\|_{\Lo[6]}^l\intd s\\&\leq \frac{2-l}{2}\intoT\|n_\epsi(\cdot,s)\|_{\Lo[\frac{6l}{6-l}]}^\frac{2l}{2-l}\intd s+\frac{l}{2}\intoT\|c_\epsi(\cdot,s)\|_{\Lo[6]}^2\intd s\quad\text{for all }\epsi\in(0,1).
\end{align*}
Due to the embedding $\W[1,2]\hookrightarrow\Lo[6]$, the bounds from Lemma \ref{lem:bounds} and \eqref{eq:strong-conv-nabc-l2-proof-eq1} entail the existence of $C_2>0$ such that
\begin{align*}
\intoTomega|n_\epsi c_\epsi|^l\leq C_2\quad\text{for all }\epsi\in(0,1),
\end{align*}
with $l=2-\frac{4}{3m}>1$. Since Lemma \ref{lem:convergences} also implies the a.e. convergence of $n_\epsi c_\epsi\to nc$ in $\Omega\times(0,\infty)$ as $\epsi=\epsi_{j_k}\searrow0$, we can employ the Vitali convergence theorem to obtain that
\begin{align*}
\intoTomega n_\epsi c_\epsi\to\intoTomega nc\quad\text{as }\epsi=\epsi_{j_k}\searrow0.
\end{align*}
Thus, making use of Lemma \ref{lem:c2-ineq}, Lemma \ref{lem:convergences} and testing the second equation of \eqref{approxprob} by $c_\epsi$ we find that
\begin{align*}
\intoTomega|\nabla c|^2&\geq-\frac{1}{2}\intomega c^2(\cdot,T)+\frac{1}{2}\intomega c_0^2-\intoTomega c^2+\intoTomega nc\\
&=\lim_{\epsi_{j_k}\searrow0}\bigg(-\frac{1}{2}\intomega c_{\epsi_{j_k}}^2(\cdot,T)+\frac{1}{2}\intomega c_0^2-\intoTomega c_{\epsi_{j_k}}^2+\intoTomega n_{\epsi_{j_k}} c_{\epsi_{j_k}}\bigg)=\lim_{\epsi_{j_k}\searrow0}\intoTomega|\nabla c_{\epsi_{j_k}}|^2.
\end{align*}
On the other hand by the lower semicontinuity of the norm in $\LSp{2}{\Omega\times(0,T)}$ with respect to weak convergence we also have
\begin{align*}
\intoTomega|\nabla c|^2\leq \liminf_{\epsi_{j_k}\searrow0}\intoTomega|\nabla c_{\epsi_{j_k}}|^2
\end{align*}
in light of \eqref{eq:conv-nab-c}. Consequently, combining the weak convergence in \eqref{eq:conv-nab-c} with the convergence of norms established above immediately implies the asserted strong convergence property.
\end{bew}

Relying on the strong convergence of $\nabla c_\epsi$ in $\LSp{2}{\Omega\times(0,T)}$ and the precompactness properties from Lemma \ref{lem:convergences}, we find that whenever $m\in(\frac{4}{3},2)$ one can check in a straightforward manner that for the choice $\Phi(s)\equiv (s+1)^{m-1}$ the supersolution property \eqref{eq:very_weak_supersol} is satisfied. Recalling that every weak solution is also a very weak solution we note that actually only $m\in(\frac43,\frac53]$ are of importance here, since for larger values Theorem \ref{theo:1} already covers the asserted existence of very weak solutions.

\begin{lemma}\label{lem:sol-prop-n}
Let $m\in(\frac{4}{3},2)$. Assume that $n_0,c_0,u_0$ comply with \eqref{IR} and denote by $n,c,u$ the limit functions provided by Lemma \ref{lem:convergences}. Moreover, $\Phi(s):=(s+1)^{m-1}$ for $s\geq0$. Then $n$ is a global $\Phi$--supersolution of \eqref{CTnonlindiff} in the sense of Definition \ref{def:weak_super_sol}.
\end{lemma}

\begin{bew}
In consideration of the regularity properties for $c$ and $u$ obtained in Lemma \ref{lem:convergences}, the fact that $m<2$, as well as \eqref{eq:conv-n+1-m-1-ae}, \eqref{eq:conv-nab-n+1-m-1}, \eqref{eq:conv-nab-n+1+eps-m-1}, \eqref{eq:conv-n-m-1-ae} and \eqref{eq:conv-nab-n-m-1} we find that all regularity requirements imposed in Definition \ref{def:weak_super_sol}, including those in \eqref{eq:very_weak_supersol_regularity}, are fulfilled and we are left with verifying that \eqref{eq:very_weak_supersol} holds. Given any nonnegative $\varphi\in C_0^\infty\big(\bomega\times[0,\infty)\big)$ with $\frac{\partial\varphi}{\partial\nu}=0$ on $\romega\times(0,\infty)$ we fix $T>0$ such that $\varphi\equiv 0$ in $\Omega\times(T,\infty)$ and test the first equation of \eqref{approxprob} by $(n_\epsi+1)^{m-2}\varphi$ to find that
\begin{align}\label{eq:sol-prop-n-proof-eq1}
&-\intoTomega(n_\epsi+1)^{m-1}\varphi_t-\intomega (n_0+1)^{m-1}\varphi(\cdot,0)\nonumber\\
=\ &m(m-1)(2-m)\!	\intoTomega \big|(n_\epsi+1)^{\frac{m-3}{2}}(n_\epsi+\epsi)^{\frac{m-1}{2}}\nabla n_\epsi\big|^2\varphi
-m\!\intoTomega (n_\epsi+\epsi)^{m-1}\big(\nabla (n_\epsi+1)^{m-1}\cdot\nabla\varphi\big)\nonumber\\
&\quad-(2-m)\intoTomega\frac{(n_\epsi+1)^{-1}n_\epsi}{(1+\epsi n_\epsi)^3}\big(\nabla (n_\epsi+1)^{m-1}\cdot\nabla c_\epsi\big)\varphi+(m-1)\intoTomega\frac{(n_\epsi+1)^{m-2}n_\epsi}{(1+\epsi n_\epsi)^3}(\nabla c_\epsi\cdot\nabla\varphi)\nonumber\\&\quad\quad
+\intoTomega (n_\epsi+1)^{m-1}(u_\epsi\cdot\nabla\varphi)\quad\text{holds for all }\epsi\in(0,1).
\end{align}
Now, since $\frac{(n_\epsi+1)^{-1}n_\epsi}{(1+\epsi n_\epsi)^3}\leq 1$ for all $\epsi\in(0,1)$ and $\frac{(n_\epsi+1)^{-1}n_\epsi}{(1+\epsi n_\epsi)^3}\to \frac{n}{n+1}$ a.e. in $\Omega\times(0,\infty)$ as $\epsi\searrow0$ we can combine the strong convergence of $\nabla c_\epsi$ in $\LSp{2}{\Omega\times(0,T)}$ obtained in Lemma \ref{lem:strong-conv-nabc-l2} with \cite[Lemma 10.4]{win15_chemorot} to discern that 
\begin{align*}
\frac{(n_\epsi+1)^{-1}n_\epsi}{(1+\epsi n_\epsi)^3}\nabla c_\epsi\to \frac{n}{n+1}\nabla c\quad\text{in }\LSp{2}{\Omega\times(0,T)}\quad\text{as }\epsi=\epsi_{j_k}\searrow0,
\end{align*}
which in sequence with \eqref{eq:conv-nab-n+1-m-1} and \eqref{eq:conv-n+1-m-1-ae} shows that
\begin{align*}
-(2-m)\intoTomega\frac{(n_\epsi+1)^{-1}n_\epsi}{(1+\epsi n_\epsi)^3}\big(\nabla(n_\epsi+1)^{m-1}\cdot\nabla c_\epsi\big)\varphi&\to-(2-m)\intoTomega\frac{n}{n+1}\big(\nabla (n+1)^{m-1}\cdot\nabla c\big)\varphi,\\
\text{and}\quad (m-1)\intoTomega\frac{(n_\epsi+1)^{m-2}n_\epsi}{(1+\epsi n_\epsi)^3}(\nabla c_\epsi\cdot\nabla\varphi)&\to(m-1)\intoTomega (n+1)^{m-2}n(\nabla c\cdot\nabla\varphi)
\end{align*}
as $\epsi=\epsi_{j_k}\searrow0$. Moreover, \eqref{eq:conv-n+1-m-1-ae}, \eqref{eq:conv-u-l2}, \eqref{eq:conv-n-m-1-ae} and \eqref{eq:conv-nab-n+1-m-1} also entail that
\begin{align*}
-\intoTomega (n_\epsi+1)^{m-1}\varphi_t&\to-\intoTomega (n+1)^{m-1}\varphi_t,\\
\intoTomega(n_\epsi+1)^{m-1}(u_\epsi\cdot\nabla\varphi)&\to \intoTomega (n+1)^{m-1}(u\cdot\nabla \varphi),\\
\text{and}\ -m\intoTomega(n_\epsi+\epsi)^{m-1}\big(\nabla(n_\epsi+1)^{m-1}\cdot\nabla\varphi\big)&\to-m\intoTomega n^{m-1}\big(\nabla (n+1)^{m-1}\cdot\nabla\varphi\big)
\end{align*}
as $\epsi=\epsi_{j_k}\!\searrow\!0$. Finally, by the lower semicontinuity of the norm in $\LSp{2}{\Omega\times(0,T)}$ with respect to weak convergence it follows from \eqref{eq:conv-nab-n+1+eps-m-1} and $m<2$ that
\begin{align*}
\liminf_{\epsi_{j_k}\searrow0}\intoTomega\big|(n_\epsi+1)^{\frac{m-3}{2}}(n_\epsi+\epsi)^{\frac{m-1}{2}}\nabla n_\epsi\big|^2\varphi\geq \intoTomega\big|(n+1)^{\frac{m-3}{2}}n^{\frac{m-1}{2}}\nabla n\big|^2\varphi,
\end{align*}
so that consolidating the statements above with \eqref{eq:sol-prop-n-proof-eq1} and the fact that $\varphi\equiv0$ in $\Omega\times(T,\infty)$ leads to
\begin{align}\label{eq:sol-prop-n-eq}
-&\intinfomega (n+1)^{m-1}\varphi_t-\intomega (n_0+1)^{m-1}\varphi(\cdot,0)\nonumber\\&\quad\geq m(m-1)(2-m)\intinfomega\big|(n+1)^{\frac{m-3}{2}}n^\frac{m-1}{2}\nabla n\big|^2\varphi-m\intinfomega n^{m-1}\big(\nabla (n+1)^{m-1}\cdot\nabla\varphi\big)\nonumber\\
&\qquad\;-(2-m)\intinfomega\frac{n}{n+1}\big(\nabla (n+1)^{m-1}\cdot\nabla c\big)\varphi+(m-1)\intinfomega (n+1)^{m-2}n(\nabla c\cdot\nabla\varphi)\\&\quad\qquad\;
+\intinfomega (n+1)^{m-1}(u\cdot\nabla\varphi)\nonumber,
\end{align}
which is equivalent to \eqref{eq:very_weak_supersol} for the choice of $\Phi(s)\equiv (s+1)^{m-1}$ and thereby concludes the proof.
\end{bew}

The proof of Theorem \ref{theo:2} is essentially finished, we only need to combine the prepared lemmas.

\begin{proof}[\textbf{Proof of Theorem \ref{theo:2}}:]
Since any global weak solution is also a global very weak solution (i.e. the $\Phi$--supersolution property is satisfied for $\Phi(s)\equiv s$), due to Theorem \ref{theo:1} we are left to treat $m\in\big(\frac{4}{3},\frac{5}{3}\big]$. For these $m$ the proof follows from an evident combination of Lemma \ref{lem:sol-prop-c-u}, Lemma \ref{lem:sol-prop-n} and the regularity properties contained Lemma \ref{lem:convergences}.
\end{proof}

\section*{Acknowledgements}
The author acknowledges support of the {\em Deutsche Forschungsgemeinschaft} in the context of the project
  {\em Analysis of chemotactic cross-diffusion in complex frameworks}. 

\footnotesize{
\setlength{\bibsep}{1pt plus 0.5ex}

}
\end{document}